\makeatletter
\let\ORIlabel\label
\let\ORIrefstepcounter\refstepcounter
\AddToHook{package/hyperref/before}{%
  \let\label\ORIlabel
  \let\refstepcounter\ORIrefstepcounter
}
\makeatother
\documentclass[onefignum,onetabnum]{siamart171218}
\usepackage{stmaryrd} 
\usepackage{xcolor}   

\usepackage[numbers]{natbib}   

\newcommand{\bn}{\mathbf{n}}
\newcommand{\bq}{\mathbf{q}}
\newcommand{\jump}[1]{\left\llbracket #1 \right\rrbracket}
\newcommand{\ave}[1]{\{ #1 \}}

\newcommand{\NN}{\mathcal{N}}
\newcommand{\LL}{\mathcal{L}}

\usepackage{amsfonts}
\usepackage{graphicx}
\usepackage{epstopdf}
\usepackage{algorithm}
\usepackage{algpseudocode} 
\newsiamremark{assumption}{Assumption}

\ifpdf
  \DeclareGraphicsExtensions{.eps,.pdf,.png,.jpg}
\else
  \DeclareGraphicsExtensions{.eps}
\fi

\newsiamremark{remark}{Remark}
\newsiamremark{hypothesis}{Hypothesis}
\crefname{hypothesis}{Hypothesis}{Hypotheses}
\newsiamthm{claim}{Claim}

\headers{DG-FEONet for Non-smooth PDEs}{K. Chawla, Y. Hong, J. Lee, and S. Lee}

\title{Discontinuous Galerkin finite element operator network for solving non-smooth PDEs\thanks{
\funding{This work was supported by the U.S. Department of Energy, Office of Science, Energy Earthshots Initiative under Award No. DE-SC0024703.}}}

\author{Kapil Chawla\thanks{Department of Mathematics, Florida State University, 1017 Academic Way, Tallahassee, FL 32306, United States of America (\email{kc25bc@fsu.edu}; \email{slee17@fsu.edu}).}
\and Youngjoon Hong\thanks{Department of Mathematical Sciences, Seoul National University, Gwanak-Ro 1, Gwanak-Gu, Seoul, Republic of Korea (\email{hongyj@snu.ac.kr}).}
\and Jae Yong Lee\thanks{Department of AI, Chung-Ang University, 84 Heukseok-ro, Dongjak-gu, Seoul, Republic of Korea (\email{jaeyong@cau.ac.kr}).}
\and Sanghyun Lee\footnotemark[2]}

\usepackage{amsopn}

\makeatletter
\newcommand*{\addFileDependency}[1]{
  \typeout{(#1)}
  \@addtofilelist{#1}
  \IfFileExists{#1}{}{\typeout{No file #1.}}
}
\makeatother

\begin{document}

\maketitle

\begin{abstract}
 We introduce Discontinuous Galerkin Finite Element Operator Network (DG--FEONet), a data-free operator learning framework that combines the strengths of the discontinuous Galerkin (DG) method with neural networks to solve parametric partial differential equations (PDEs) with discontinuous coefficients and non-smooth solutions. Unlike traditional operator learning models such as DeepONet and Fourier Neural Operator, which require large paired datasets and often struggle near sharp features, our approach minimizes the residual of a DG-based weak formulation using the Symmetric Interior Penalty Galerkin (SIPG) scheme. DG-FEONet predicts element-wise solution coefficients via a neural network, enabling data-free training without the need for precomputed input-output pairs. We provide theoretical justification through convergence analysis and validate the model's performance on a series of one- and two-dimensional PDE problems, demonstrating accurate recovery of discontinuities, strong generalization across parameter space, and reliable convergence rates. Our results highlight the potential of combining local discretization schemes with machine learning to achieve robust, singularity-aware operator approximation in challenging PDE settings.
\end{abstract}

\begin{keywords}
  discontinuous Galerkin methods , finite element methods  , physics-informed operator learning  ,  discontinuity-aware learning
\end{keywords}

\begin{AMS}
 65M60, 65N30, 68T20, 68U07
\end{AMS}

\section{Introduction}
Partial differential equations (PDEs) form the mathematical foundation for modeling complex physical systems, from fluid flow and heat transfer to electromagnetics and elasticity~\cite{evans2022partial}. While analytical solutions exist for some simple cases, real-world PDE problems typically involve irregular domains, nonlinearities, or parametric uncertainty, making numerical methods essential for practical computation~\cite{brenner2008mathematical}. Among these, the finite element method (FEM) stands out for its flexibility and accuracy, particularly in handling complex geometries and boundary conditions.

Despite their reliability, conventional numerical solvers can be computationally expensive when repeatedly applied to parametric PDEs-for example, in design optimization, uncertainty quantification, or control applications. This has led to a surge of interest in machine learning-based methods that can approximate the solution operator itself. Early approaches, including Physics-Informed Neural Networks (PINNs)~\cite{MR3881695}, offer data-free learning by embedding PDE residuals into the loss function. However, PINNs often struggle with sharp gradients, boundary layers, and discontinuous coefficients~\cite{liu2024discontinuity,mao2023physics,SP01}. Moreover, they are typically trained for one specific problem instance, limiting their generalizability.

Recently, the concept of neural operators has gained attention as a way to directly learn solution operators that change with varying parameters in parametric PDEs~\cite{LU2022114778,kovachki2021neural}. Neural operator frameworks such as DeepONet~\cite{lu2021learning} and the Fourier Neural Operator (FNO)~\cite{li2021fourier} aim to learn mappings between infinite-dimensional function spaces, enabling efficient prediction of solutions without the need to re-solve the PDE for each new parameter setting. Beyond the DeepONet and the FNO, various neural operators have been proposed, including graph-based~\cite{brandstetter2022message,pfaff2021learning}, wavelet-based~\cite{gupta2021multiwaveletbased}, and transformer-inspired models~\cite{cao2021choose,wang2025cvit,pmlr-v235-hao24d}, aiming to address irregular domains and improve generalization across parameter spaces. However, they generally rely on supervised training using large datasets of input-output solution pairs, which can be expensive or impractical to generate, especially for high-dimensional or nonlinear problems.

To overcome these challenges, recent studies have explored incorporating physics-based losses, inspired by the approach of PINNs, into the operator learning framework~\cite{Goswami2023}. For example, \cite{wang2021learning} and \cite{PINO} introduced physics-informed residual losses into DeepONet and FNO architectures, respectively, enabling the networks to leverage governing PDE structures during training. Building on this line of work, the Finite Element Operator Network (FEONet)~\cite{MR4888707,hong2024error} was recently proposed as a physics-aware, data-free alternative. FEONet integrates the continuous Galerkin (CG) formulation of the FEM into the neural network training loop by minimizing the weak form residual instead of the solution error. By predicting FEM coefficients and reconstructing the solution via basis functions, FEONet eliminates the need for labeled data and inherits essential FEM properties such as mesh flexibility and exact enforcement of boundary conditions.

However, the CG formulation enforces continuity across element boundaries and does not inherently satisfy local conservation properties. This makes it ill-suited for problems involving discontinuities in the solution or coefficients, such as shock waves, fractured media, or multi-phase flow. In such settings, CG methods can introduce spurious oscillations and fail to capture localized features accurately~\cite{COCKBURN1998199}. In this work, we propose DG-FEONet, a novel extension of FEONet based on the discontinuous Galerkin (DG) method using the Symmetric Interior Penalty Galerkin (SIPG) formulation. The DG framework relaxes the inter-element continuity requirement and enforces consistency through numerical fluxes and penalty terms. This makes DG methods particularly effective for solving PDEs with discontinuities, non-smooth solutions, or heterogeneous media, while preserving conservation laws locally.

Our approach predicts the local DG coefficients element-wise using a neural network and defines the loss based on the SIPG residual. Like the original FEONet, DG-FEONet requires no precomputed solutions and is trained in an data-free manner. However, by leveraging the DG framework, it significantly improves accuracy and robustness in problems where CG-based models struggle. Through a series of numerical experiments, we demonstrate that DG-FEONet accurately captures sharp features and discontinuities in parametric PDEs, and offers a reliable, data-free approach for operator learning in challenging regimes.

The remainder of this paper is organized as follows. Section~\ref{sec:related_work} provides an overview of related work. In Section~\ref{sec:model}, we present the formulation of DG-FEONet, including a brief overview of the DG-FEM and the construction of the weak-form loss based on the Symmetric Interior Penalty Galerkin (SIPG) scheme. Section~\ref{sec:analysis} provides a theoretical convergence analysis of the proposed framework. In Section~\ref{sec:experiments}, we validate DG-FEONet through a series of numerical experiments involving both one-dimensional and two-dimensional parametric PDEs with discontinuous coefficients. We also evaluate the convergence behavior and demonstrate the robustness of our method across varying mesh resolutions. Finally, Section~\ref{sec:conclusion} concludes the paper and outlines directions for future work.

\section{Related works}\label{sec:related_work}
The FEM has long been established as a cornerstone for solving PDEs numerically~\cite{FEM1, FEM2}. The CG formulation is the classical FEM approach, where the solution is approximated using globally continuous basis functions. It is particularly effective for problems with smooth solutions. However, its enforcement of inter-element continuity limits its applicability to problems featuring sharp gradients or discontinuities. The DG method~\cite{shu2009discontinuous,riviere2008discontinuous} overcomes these limitations by allowing the solution to be discontinuous across element boundaries. DG methods have been widely adopted for conservation laws, and multiphase flow, owing to their local conservation properties, flexibility in mesh handling, and compatibility with high-order approximations. Among various DG variants, the Symmetric Interior Penalty Galerkin (SIPG)method
\cite{muller2017symmetric,riviere2008discontinuous} is a popular choice for elliptic and parabolic problems due to its stability and ease of implementation.

Recent advances in deep learning have positioned neural networks as powerful tools for tackling singularly perturbed PDEs. Physics-informed machine learning has also been explored for these challenges~\cite{SP01, SP02,chen2025learn}, but typically lacks the scalability needed for more general applications. To better handle discontinuities, recent works have proposed PINN-based methods specifically targeting discontinuous solutions~\cite{liu2024discontinuity, wang2025discontinuous}. In the context of neural operators, component FNO~\cite{SP03} addresses singular perturbations by incorporating asymptotic expansions into its architecture, offering a promising direction for improving the robustness and accuracy of neural network-based solvers.

Recent developments have explored hybrid methodologies that integrate classical numerical methods with neural networks to leverage the strengths of both paradigms. These approaches aim to improve the efficiency, generalizability, and physical consistency of PDE solvers. For instance, neural networks have been used to learn data-driven discretization schemes~\cite{bar2019learning}, effectively replacing traditional finite difference operators. FEONet~\cite{MR4888707,hong2024error} combines deep learning with the CG-FEM by directly learning expansion coefficients of the solution in the FEM basis, eliminating the need for labeled data. Similarly, spectral-based approaches such as Spectral Coefficient Learning via Operator Network (SCLON)\cite{CHOI2024116678} and Unsupervised Legendre Galerkin Network (ULGNet)\cite{10043713} learn solution coefficients in spectral or Galerkin bases, enabling physics-aware data-free training. More recently, the Neural-Operator Element Method (NOEM)\cite{ouyang2025neural} was introduced to reduce meshing costs in multiscale problems by embedding neural operators as surrogate elements within a variational FEM framework. Recently, \cite{chen2025dgnn, shen2025hybrid, franck2024approximately} have proposed DG-based neural frameworks, where DG formulations are integrated with neural networks either as trial/test spaces, hybrid solvers, or enriched bases.

\section{Model Description}{\label{sec:model}}

We consider the problem of approximating the solution operator of parametric PDEs using a DG-FEM enhanced with operator learning. In particular, we consider the following governing equation:
\begin{equation} \label{eq:pde}
- \nabla \cdot \left( \varepsilon(x) \nabla u \right) + F(u) = f(x), \quad x \in \Omega,
\end{equation}
in the  domain $\Omega\subset\mathbb{R}^d$ ($d=1,2,3$), where the scalar solution $u(x)\colon \Omega \to \mathbb{R}$ satisfies appropriate boundary conditions on \(\partial\Omega\), where the boundary is decomposed into its Dirichlet and Neumann parts:
$
\partial\Omega = \partial\Omega_D \cup \partial\Omega_N, \text{ and }
\partial\Omega_D \cap \partial\Omega_N = \emptyset.$ As a representative example, we consider
$$
 F(u)  = b(x) \nabla u(x) + c(x) u(x),
$$
which is often referred to as the convection–diffusion–reaction equation. Here, $\varepsilon(x)$, $b(x)$, and $c(x)$ are spatially varying coefficients, $f(x)$ is the prescribed source function, and these quantities may be discontinuous across $\Omega$. Our goal is to learn the underlying solution operator
\begin{equation}
    \mathcal{G} : (\varepsilon, b, c, f) \;\mapsto\; u,
\end{equation}
which maps a collection of spatially varying coefficients and source terms to the corresponding solution $u(x)$ of the PDE~\eqref{eq:pde}. Any of the fields $\varepsilon(x)$, $b(x)$, $c(x)$, or $f(x)$ may serve as the parametric input, depending on the problem setting, and DG-FEONet aims to approximate $\mathcal{G}$ for the chosen input configuration.
\subsection{Discontinuous Galerkin Finite Element Methods}

Unlike classical CG-FEM, the DG method~\cite{wheeler1978elliptic,riviere1999improved,riviere2008discontinuous,cockburn2012discontinuous} allows both the trial and test spaces to be discontinuous across element boundaries. Let $\mathcal{T}_h$ be a shape-regular partition of the domain $\Omega$ into elements $K$.
The set of all edges in $\mathcal{T}_{h}$ is denoted by $\varepsilon_{h} = \varepsilon^{\text{int}}_{h}\cup\varepsilon^D_{h} \cup\varepsilon^N_{h}$, which is decomposed into three disjoint subsets. Here, $\varepsilon^{\text{int}}_{h}$ is the set of all interior edges, $\varepsilon^{D}_h$ and $\varepsilon^N_h$ are the set of all boundary edges where Dirichlet conditions or Neumann conditions are imposed, respectively.  For convenience, we also denote $\varepsilon^{I}_{h} := \varepsilon^{\text{int}}_{h}\cup \varepsilon^{D}_{h}$ as the union of the interior and Dirichlet boundary edges.

For each edge $e\in\varepsilon_h$, let $K^{+}$ and $K^{-}$ be the neighboring elements such that $e=\partial K^{+}\cap\partial K^{-}$. Let $\bn^{+}$ and $\bn^{-}$ be the unit outward normal vector to $\partial K^{+}$ and $\partial K^{-}$ respectively.
For a given vector function $\mathbf{q}$ and scalar function $v$, we define the weighted average operator $\{\cdot\}$ and the jump operator $\left\llbracket \cdot\right\rrbracket $ by
$$
\{\bq\} = \frac{\bq|_{K^{+}} + \bq|_{K^{-}}}{2}\;\text{ on } e, \quad \quad 
\left\llbracket v\right\rrbracket =v|_{K^{+}}\bn^{+}+v|_{K^{-}}\bn^{-}\;\text{on } e.
$$
If $e \in \partial \Omega$, then $e$ belongs to only one element $K$, and  we define
$$
\jump{v} = \ave{v}  =  (v \vert_{K})\vert_e.
$$

For the DG spatial discretization, we define the discrete solution space $V_h \subset L^2(\Omega)$ as a broken polynomial space:
\begin{equation}\label{DG_space}
V_h := \{ v \in L^2(\Omega) \mid v|_K \in \mathbb{P}_r(K), \ \forall K \in \mathcal{T}_h \},
\end{equation}
where $\mathbb{P}_r(K)$ denotes the space of polynomials of degree at most $r$ on element $K$.

The discrete solution $u_h \in V_h$ is determined by the Symmetric Interior Penalty Galerkin (SIPG) formulation~\cite{sun2005discontinuous}: find $u_h \in V_h$ such that
\begin{equation} \label{eq:sipg}
a(u_h, v_h) = \ell(v_h) \quad \forall v_h \in V_h,
\end{equation}
where the bilinear form $a(\cdot, \cdot)$ is given as

\begin{align}
a(u, v) 
&= \sum_{K \in \mathcal{T}_h} \int_K 
\varepsilon \nabla u \cdot \nabla v 
+ (b \cdot \nabla u) v 
+ c u v \, dx \nonumber \\
&\quad - \sum_{e \in \mathcal{E}_h^\text{int}} \int_e 
\left\{ \! \! \{ \varepsilon \nabla u \} \! \! \right\} \cdot \llbracket v \rrbracket 
+ \left\{ \! \! \{ \varepsilon \nabla v \} \! \! \right\} \cdot \llbracket u \rrbracket \, ds \nonumber \\
&\quad + \sum_{e \in \mathcal{E}_h^\text{int}} \int_e 
\frac{\sigma}{h_e} \llbracket u \rrbracket \cdot \llbracket v \rrbracket \, ds \nonumber  - \sum_{e \in \mathcal{E}_h^\text{int}} \int_e 
\widehat{b u} \cdot \llbracket v \rrbracket \, ds 
- \sum_{e \in \mathcal{E}_h^{\partial}} \int_e 
\widehat{b u} \cdot v \, ds \nonumber \\
&\quad - \int_{\partial \Omega_D} \varepsilon \nabla u \cdot \mathbf{n} \, v \, ds
- \int_{\partial \Omega_D} \varepsilon \nabla v \cdot \mathbf{n} \, u \, ds
+ \int_{\partial \Omega_D} \frac{\sigma}{h} u v \, ds.
\end{align}

Here, $ \llbracket \cdot \rrbracket $ denotes the jump across interior faces, $ \{\!\!\{ \cdot \}\!\!\} $ the average, $ h_e $ the local mesh size, and $ \sigma > 0 $ the penalty parameter. The numerical flux $ \widehat{b u} $ for the convection term is defined as $$\widehat{b u} = \frac{1}{2}(b^+ + b^-) \{ u \} + \frac{1}{2} \left| b^+ + b^- \right| \llbracket u \rrbracket,$$
The first term provides consistency, while the second introduces numerical dissipation for stability. 

The linear form is defined by:

\begin{equation}\label{eq:load}
\begin{aligned}
\ell(v)
&= \sum_{K \in \mathcal{T}_h} \int_K f(x)\,v(x)\,dx
 - \sum_{e \in \mathcal{E}_h^{D}} \int_e \bigl(\varepsilon \nabla v \cdot \bn \bigr)\,u_D\,ds \\
&\quad + \sum_{e \in \mathcal{E}_h^{D}} \int_e \frac{\sigma}{h_e}\,u_D\,v\,ds
 + \sum_{e \in \mathcal{E}_h^{N}} \int_e u_N\,v\,ds .
\end{aligned}
\end{equation}

where $u_D$ is given Dirichelt boundary condiition on $\partial \Omega_D$ and $u_N$ is Neumann bonudary condition on $\partial \Omega_N$.

To represent the DG solution in a finite-dimensional form, we express \( u_h \in V_h \) as a linear combination of locally supported basis functions \( \{ \phi_k(x) \}_{k=1}^{N(h)} \):
$$
u_h(x) = \sum_{k=1}^{N(h)} \alpha_k(\omega) \phi_k(x),
$$
where \( \alpha_k(\omega) \in \mathbb{R} \) are the expansion coefficients that depend on the problem data \( \omega \), such as the diffusion coefficient \( \varepsilon(x) \), the source term \( f(x) \), or other parametric inputs. These coefficients are typically computed by solving the linear system arising from the DG variational formulation.

Solving the linear systems arising from DG discretizations is computationally demanding and often nontrivial due to the increased degrees of freedom and the block structure of the resulting matrices~\cite{dobrev2006two,karakashian2018two}. The FEONet framework offers a promising alternative by learning the mapping from problem parameters to solution coefficients directly, thereby circumventing the need for explicit linear system assembly and solution.

\subsection{Discontinuous Galerkin Finite Element Operator Network}

Solving a parametric PDE using traditional DG methods requires assembling and solving a linear system for each new sample of the input parameters (e.g., coefficients, source terms). Instead of repeatedly solving this system, DG-FEONet aims to learn a direct map from problem inputs to the solution coefficients.

Let $\{ \phi_k \}_{k=1}^{N(h)}$ denote a local DG basis for the discrete space $V_h$. Then the approximate solution is represented as:
$$
u_h(x; \omega) = \sum_{k=1}^N \alpha_k(\omega) \phi_k(x),
$$
where $\omega \in \Omega_\text{param}$ denotes a sample of parametric inputs. The goal is to learn the map $\omega \mapsto \alpha(\omega)$ using a neural network $\mathcal{N}_\theta$:
$$
\alpha(\omega) \approx \mathcal{N}_\theta(\omega).
$$

This idea draws inspiration from PINNs, which incorporate physical laws into the learning process by minimizing PDE residuals. However, unlike classical PINNs that approximate the solution directly at collocation points, DG-FEONet operates in a \textit{finite-dimensional function space}. The neural network predicts the coefficients of a DG basis expansion, effectively learning a parametric solution operator. This aligns DG-FEONet more closely with operator learning frameworks such as the FEONet, DeepONet and FNO, but with the distinct advantage of using DG discretization for better handling of discontinuities.

Another key distinction is that DG-FEONet is trained in an data-free manner, without requiring access to labeled solutions. Instead, we define a physics-informed loss based on the residual of the DG variational formulation:
\begin{equation} \label{eq:loss}
\mathcal{L}(\theta) = \frac{1}{M} \sum_{m=1}^M \sum_{i=1}^N \left| a(u_h^{(m)}, \phi_i) - \ell^{(m)}(\phi_i) \right|^2,
\end{equation}
where $u_h^{(m)}(x) = \sum_{k=1}^N \mathcal{N}_\theta(\omega_m)_k \phi_k(x)$ is the predicted solution for the $m$-th input. The residual is evaluated against test functions $\phi_i$ from the DG basis.

\begin{figure}[t]
    \centering
    \includegraphics[width=\textwidth]{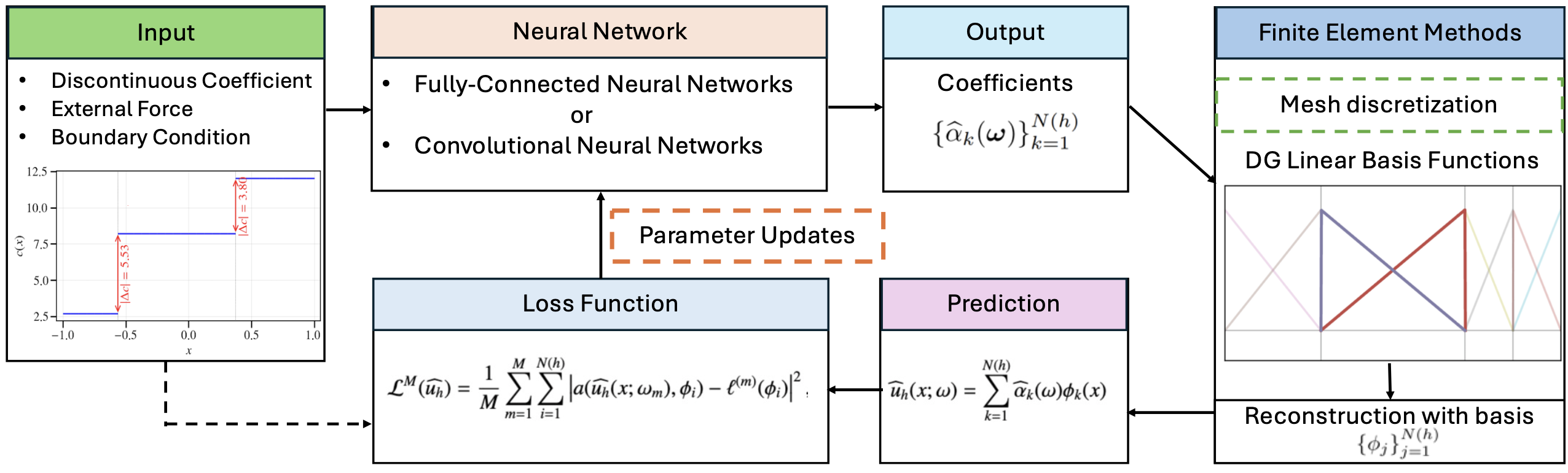}
    \caption{ {Schematic overview of the DG-FEONet framework. The PDE inputs, which we emphasize that these can be discontinuous, are passed to a neural network that predicts the coefficients for DG-FEM solutions. These coefficients are used to reconstruct the DG solution using basis functions defined over a triangulated domain. The network is trained in an data-free manner using the residual of the DG weak formulation as the loss function. }}
    \label{fig:dgfeonet-architecture}
\end{figure}

This formulation enables \textit{data-free operator learning} directly from the physics of the PDE, leveraging the structure of DG methods. The locality and flexibility of DG bases make this framework especially suitable for problems with discontinuous coefficients or solutions, where CG-based models cannot handle discontinuous solutions. Figure~\ref{fig:dgfeonet-architecture} provides a schematic overview of the DG-FEONet framework.

\subsection{Loss Formulation and Training Objective}\label{loss_section}

The DG-FEONet framework is trained by minimizing the residual of the weak formulation associated with the SIPG discretization. Given a neural network prediction of the DG solution in the form
$$
u_h(x; \omega) = \sum_{k=1}^{N} \hat{\alpha}_k(\omega) \phi_k(x),
$$
where $ \hat{\alpha}_k(\omega) $ are the network-predicted coefficients and \( \{ \phi_k \}_{k=1}^N \) is a local DG basis, we define the training objective by enforcing the discrete variational form:
$$
a(u_h^{(m)}, \phi_i) \approx \ell^{(m)}(\phi_i), \quad \forall i = 1, \dots, N,
$$
for each sample $ \omega_m \in \Omega_{\text{param}} $, $ m = 1, \dots, M $. Here, $ a(\cdot, \cdot) $ is the SIPG bilinear form that includes volume integrals, numerical flux terms across interior faces, and penalty terms, as previously defined in Equation~\eqref{eq:sipg}.

We define the residual for the \( m \)-th input sample and the \( i \)-th test function as
\[
\mathcal{R}_i^{(m)} := a(u_h^{(m)}, \phi_i) - \ell^{(m)}(\phi_i),
\]
and formulate the training loss as the average squared residual over all input samples and test functions:
\begin{equation}\label{loss_equation}
\mathcal{L}_{\text{DG}}(\theta) = \frac{1}{M} \sum_{m=1}^{M} \sum_{i=1}^{N} \left| \mathcal{R}_i^{(m)} \right|^2.
\end{equation}

This loss function measures how well the neural network prediction satisfies the discrete variational statement for each sample in an data-free manner. No labeled solutions \( u(x; \omega) \) are needed, making this approach fully physics-informed. By minimizing the residuals across the DG test space, the network learns a mapping from input functions (e.g., coefficients or source terms) to solution coefficients that approximate the solution operator of the PDE.

To ensure computational efficiency, we pre-assemble the SIPG matrices corresponding to the stiffness, convection, and interface flux terms. During training, the neural network outputs the predicted coefficients \( \hat{\alpha}(\omega_m) \), which are then inserted into the assembled system to compute residuals and their gradients for backpropagation.

The overall training workflow is summarized in Algorithm~\ref{alg:dgfeonet}, which outlines the key steps of the DG-FEONet optimization procedure.

\begin{algorithm}[t]
\caption{DG-FEONet Training Procedure (SIPG Formulation)}
\label{alg:dgfeonet}
\begin{algorithmic}[1]
\Require Training samples $\{f(x;\omega_m)\}_{m=1}^M$
\Ensure Predicted coefficients $\widehat{\alpha}(\omega_m) = [\widehat{\alpha}_k(\omega_m)]_{k=1}^{N(h)}$ and predicted solution
$\widehat{u}_h(x;\omega_m)=\sum_{k=1}^{N(h)} \widehat{\alpha}_k(\omega_m)\varphi_k(x)$
\State \textbf{Pre-compute:}
\State Assemble SIPG stiffness matrix $A^{(s)}$ (penalty/jump terms).
\State Assemble convection matrix $A^{(c)}$ (numerical fluxes and boundary terms).
\State Assemble load vectors $F^{(m)}$ for $m=1,\dots,M$ using $f(x;\omega_m)$.
\For{$n = 1,2,\dots,N_{\text{epochs}}$}
  \State $\mathcal{L}_{\mathrm{DG}} \gets 0$
  \For{$m = 1,2,\dots,M$}
    \State $\widehat{\alpha}(\omega_m) \gets \texttt{DG-FEONet}(f(x;\omega_m))$
    \State $\mathcal{R}^{(m)} \gets (A^{(s)} + A^{(c)})\widehat{\alpha}(\omega_m) - F^{(m)}$
    \State $\mathcal{L}_{\mathrm{DG}} \gets \mathcal{L}_{\mathrm{DG}} + \|\mathcal{R}^{(m)}\|^2$
  \EndFor
  \State Update DG-FEONet parameters to minimize $\mathcal{L}_{\mathrm{DG}}$
\EndFor
\end{algorithmic}
\end{algorithm}

\section{Convergence Analysis of DG-FEONet}\label{sec:analysis}Throughout this section, we use the DG setting introduced in
Section \ref{sec:model}. In particular, $\Omega \subset \mathbb{R}^d$ denotes the physical domain,
$\Omega_{\text{param}}$ the parameter space equipped with probability measure $\mathbb{P}_{\Omega_{\text{param}}}$,
$V_h$ the broken DG space defined in \eqref{DG_space}, and
$a(\cdot,\cdot)$ the SIPG bilinear form.

\begin{definition}[DG energy norm]\label{def:DGnorm_main}
For $v_h \in V_h$, define
\begin{equation}\label{eq:DGnorm_main}
\begin{aligned}
\|v_h\|_{\mathrm{DG}}^2
:=\;&
\sum_{K \in \mathcal{T}_h}
\Bigl(
  \,\|\sqrt{\varepsilon}\nabla v_h\|_{L^2(K)}^2
  + \|v_h\|_{L^2(K)}^2
\Bigr) \;+
\sum_{e \in \mathcal{E}_h^{\mathrm{int}}}
\frac{\sigma}{h_e}\,
\|\llbracket v_h \rrbracket\|_{L^2(e)}^2 \\
&\;+
\sum_{e \subset \partial\Omega}
\frac{\sigma}{h_e}\,
\|v_h\|_{L^2(e)}^2 .
\end{aligned}
\end{equation}
\end{definition}


\begin{assumption}[Coercivity]\label{ass:coercivity_main}
Assume that the coefficients in the convection-diffusion-reaction operator introduced in Section~\ref{sec:model} satisfy the standard positivity condition
\begin{equation}\label{eq:positivity_condition}
c(x) - \tfrac12 \nabla\cdot b(x) \ge c_0 \ge 0
\qquad \text{a.e. in } \Omega,
\end{equation}
and that the SIPG penalty parameter $\sigma>0$ is chosen sufficiently large.
Then there exists a constant $C_{\mathrm{coer}}>0$, independent of $h$, such that
\begin{equation}\label{eq:coercivity_main}
a(v_h,v_h)\ge C_{\mathrm{coer}}\|v_h\|_{DG}^2
\qquad \forall v_h\in V_h.
\end{equation}
\end{assumption}

The coercivity of the SIPG diffusion term with sufficiently large penalty
parameter is well established; see, e.g., \cite{EpshteynRiviere2007}.
For convection-diffusion-reaction problems, stability of the DG bilinear
form under the positivity condition \eqref{eq:positivity_condition} follows
from standard DG theory; see, e.g., \cite{Riviere2008}.


\subsection{Coefficient representation}\label{subsec:coef_rep}
Let $\{\phi_i\}_{i=0}^{N_h}$ be a basis of $V_h$. We define the DG operator matrix
$A\in\mathbb{R}^{(N_h+1)\times(N_h+1)}$ associated with $a(\cdot,\cdot)$ by
\begin{equation}\label{eq:A_def_main}
A_{ij} := a(\phi_j,\phi_i).
\end{equation}
Under Assumption~\ref{ass:coercivity_main}, $A$ is invertible; see Appendix \ref{app:A_invertible}.

\subsection{Approximation error}

\begin{lemma} \cite{HornJohnson2012}\label{matrix_thm_gen}
Let \( A \in \mathbb{R}^{(N_h + 1) \times (N_h + 1)} \) be invertible (not necessarily symmetric or positive definite). Denote its minimal and maximal singular values by \(\sigma_{\min}(A)\) and \(\sigma_{\max}(A)\), respectively. Then, for any \(\mathbf{x} \in \mathbb{R}^{N_h + 1}\),
$$
\sigma_{\min}(A)\,\|\mathbf{x}\|_2 \;\le\; \|A\mathbf{x}\|_2 \;\le\; \sigma_{\max}(A)\,\|\mathbf{x}\|_2.
$$
Equivalently, letting \(S = A^\top A\) (which is symmetric positive definite), we obtain
$$
\lambda_{\min}(S) = \sigma_{\min}(A)^2, \; \lambda_{\max}(S) = \sigma_{\max}(A)^2, \; 
$$
\text{ and } \;
$$
\lambda_{\min}(S)\,\|\mathbf{x}\|_2^2 \;\le\; \|A\mathbf{x}\|_2^2 \;\le\; \lambda_{\max}(S)\,\|\mathbf{x}\|_2^2.
$$
\end{lemma}

Before stating the approximation result, we introduce the neural network
approximation in coefficient space. Let $\mathcal{N}_n$ denote the class of
neural networks with capacity parameter $n$. Let $F\in\mathbb{R}^{N_h+1}$ be
the DG load vector with entries
$$
F_i := \ell(\phi_i), \qquad i=0,\dots,N_h,
$$
where $\{\phi_i\}_{i=0}^{N_h}$ is a basis of $V_h$ and $\ell(\cdot)$ is the
linear functional defined in \eqref{eq:load}. We define the loss functional
$$
\mathcal{L}(\alpha) := \|A\alpha - F\|_2^2,
$$
and the population minimizer
\begin{equation}\label{def_alpha}
\widehat{\alpha}(n)
:= \arg\min_{\alpha \in \mathcal{N}_n} \mathcal{L}(\alpha).
\end{equation}
With this notation in place, we now establish convergence of the approximation
error in coefficient space.

\begin{theorem}\cite{MR4888707}\label{approx_main_thm}
Let \( \alpha^* \) denote the true coefficient function, and let \( \widehat{\alpha}(n) \) be the approximation obtained from the DG-FEONet model trained with \( n \) capacity parameter independent of the sample size $M$. Then the following convergence holds:
\begin{equation}\label{app_t_est}
    \| \alpha^* - \widehat{\alpha}(n) \|_{L^2(\Omega_{\text{param}})} \to 0 \quad \text{as } n \to \infty.
\end{equation}
\end{theorem}

\subsection{Generalization Error}
We use the concept of Rademacher complexity, which serves as a quantitative measure of the capacity of a function class to fit purely random noise \cite{Rade_1, Rade_2}. 

\begin{definition}
Let $\{X_i\}_{i=1}^M$ be a collection of i.i.d.~random variables. The Rademacher complexity of a function class $\mathcal{G}$ is given by
\[
    R_M(\mathcal{G}) \;=\; 
    \mathbb{E}_{\{X_i,\varepsilon_i\}_{i=1}^M}\!\left[\;
    \sup_{f\in\mathcal{G}}
    \bigg|\frac{1}{M}\sum_{i=1}^M \varepsilon_i f(X_i)\bigg|
    \;\right],
\]
where $\{\varepsilon_i\}_{i=1}^M$ denotes a sequence of i.i.d.~Rademacher random variables, that is, $\mathbb{P}(\varepsilon_i = 1)=\mathbb{P}(\varepsilon_i=-1)=\tfrac{1}{2}$ for each $i=1,\dots,M$.
\end{definition}

\begin{theorem}\label{rade_thm}
Let $\mathcal{G}$ be a $b$-uniformly bounded class of functions. Then, for any $\delta>0$, one has
\[
\sup_{f\in\mathcal{G}}
\bigg|\frac{1}{M}\sum_{i=1}^M f(X_i) - \mathbb{E}[f(X)]\bigg|
    \;\leq\; 2R_M(\mathcal{G})+\delta,
\]
with probability at least 
$
1-\exp\!\Big(-\tfrac{M\delta^2}{2b^2}\Big).
$
\end{theorem}

\begin{proof}
The proof is a direct consequence of Theorem~4.10 in \cite{Rade_2}.
\end{proof}

Next, we introduce the following function class:
\begin{equation}\label{fct_class}
    \mathcal{G}_n := \{\, |A\alpha - F|^2 : \alpha \in \NN_n \,\}.
\end{equation}
From Lemma~\ref{matrix_thm_gen}, there exists a constant $C>0$ such that
\begin{equation}\label{eq:triangle_bound}
\begin{aligned}
\|A\alpha - F\|_{L^{\infty}(\Omega_{\text{param}})}
&\le \|A\alpha\|_{L^{\infty}(\Omega_{\text{param}})}
   + \|F\|_{L^{\infty}(\Omega_{\text{param}})} \\
&\le C\Bigl(
     \|\alpha\|_{L^{\infty}(\Omega_{\text{param}})}
     + \|f\|_{C(\Omega_{\text{param}};L^1(\Omega))}
     \Bigr).
\end{aligned}
\end{equation}

The subsequent lemma is an immediate consequence of Theorem~\ref{rade_thm} in our setting.

\begin{lemma}\label{rade_main}
Let $\{\omega_m\}_{m=1}^M$ be i.i.d.~samples drawn from the distribution $\mathbb{P}_{\Omega_{\text{param}}}$. Then, for any $\delta > 0$, with probability at least $1 - 2 \exp\!\big(-\tfrac{M\delta^2}{32 \tilde{b}^2}\big)$, one has
\begin{equation}\label{loss_diff}
    \sup_{\alpha \in \NN_n} 
    \Big| \LL^M(\alpha) - \LL(\alpha) \Big|
    \;\leq\; 2R_M(\mathcal{G}_n) + \tfrac{\delta}{2}.
\end{equation}
\end{lemma}

We now employ Lemma~\ref{rade_main} to establish a convergence result for the generalization error. Throughout, we assume that the Rademacher complexity of $\mathcal{G}_n$ vanishes as $M \to \infty$, which is valid in a variety of settings \cite{Rade_1, Rade_2}.

\begin{theorem}\cite{MR4888707}\label{gen_conv_thm}
For any $n \in \mathbb{N}$, suppose that $\lim_{M\to\infty} R_M(\mathcal{G}_n) = 0$. Then, with probability one, it holds that
\[
    \lim_{n \to \infty} \, \lim_{M \to \infty} 
    \|\widehat{\alpha}(n,M) - \widehat{\alpha}(n)\|_{L^2(\Omega_{\text{param}})} = 0.
\]
\end{theorem}
\begin{lemma}\cite{HornJohnson2012}\label{lem:SPD_quadratic_bound}
Let $B \in \mathbb{R}^{(N_h+1)\times(N_h+1)}$ be a symmetric positive definite matrix, 
and denote its smallest and largest eigenvalues by $\lambda_{\min}(B)$ and $\lambda_{\max}(B)$, respectively. Then, for every vector $x \in \mathbb{R}^{N_h+1}$,
\begin{equation}\label{eq:SPD_quadratic_bound}
    \lambda_{\min}(B)\,\|x\|_2^2 
    \;\le\; x^\top B x 
    \;\le\; \lambda_{\max}(B)\,\|x\|_2^2,
\end{equation}
where $\|\cdot\|_2$ denotes the Euclidean norm on $\mathbb{R}^{N_h+1}$.
\end{lemma}

\subsection{Convergence of DG-FEONet}
By combining Theorem~\ref{approx_main_thm} with Theorem~\ref{gen_conv_thm}, we immediately obtain
\begin{equation}\label{conv_2}
    \lim_{n \to \infty} \,
    \lim_{M \to \infty} \|\alpha^* - \widehat{\alpha}(n,M)\|_{L^2(\Omega_{\text{param}})} = 0.
\end{equation}

\begin{theorem}\label{main_thm_whole}
Let $n \in \mathbb{N}$ and let $h>0$, $\varepsilon>0$ be fixed.
Assume that, for each $n$, the Rademacher complexity of
$$
    \widetilde{\mathcal{G}}_n := \{\, |A\alpha - F|^2 : \alpha \in \NN_n \,\}
$$
satisfies $R_M(\widetilde{\mathcal{G}}_n) \to 0$ as $M \to \infty$, so that the coefficient space generalization result in Theorem~\ref{gen_conv_thm} holds.
Then, with probability one,
\begin{equation}\label{main_conv_whole}
    \lim_{n \to \infty} \, \lim_{M \to \infty}
    \big\| u^{\rm DG}_{\varepsilon,h} - \widehat{u}^{\rm DG}_{\varepsilon,h,n,M}
    \big\|_{L^2(\Omega_{\text{param}};L^2(\Omega))} = 0.
\end{equation}
\end{theorem}

\begin{proof}
    Let $\{\phi_i\}_{i=0}^{N(h)}$ be a basis of $V_h$.For each parameter $\omega\in\Omega_{\text{param}}$, we can write the DG solution and its DG-FEONet approximation in coefficient form as
$$
  u^{\mathrm{DG}}_{\varepsilon,h}(\omega,x)
   \;=\; \sum_{i=0}^{N(h)} \alpha_i^*(\omega)\,\phi_i(x),
  \qquad
  \widehat{u}^{\mathrm{DG}}_{\varepsilon,h,n,M}(\omega,x)
   \;=\; \sum_{i=0}^{N(h)} \widehat{\alpha}_i(n,M;\omega)\,\phi_i(x).
$$
Define the coefficient difference
$$
  \Delta\alpha(\omega)
  \;:=\; \alpha^*(\omega)-\widehat{\alpha}(n,M;\omega)
  \in \mathbb{R}^{N(h)+1},
$$
and the corresponding DG function difference
$$
  \Delta u^{\mathrm{DG}}_{\varepsilon,h}(\omega,x)
  \;:=\; u^{\mathrm{DG}}_{\varepsilon,h}(\omega,x)
       - \widehat{u}^{\mathrm{DG}}_{\varepsilon,h,n,M}(\omega,x)
  \;=\; \sum_{i=0}^{N(h)} \Delta\alpha_i(\omega)\,\phi_i(x).
$$

For each fixed $\omega\in\Omega_{\text{param}}$, we have
\begin{align*}
  \|\Delta u^{\mathrm{DG}}_{\varepsilon,h}(\omega,\cdot)\|_{L^2(\mathcal T_h)}^2
  &:= \sum_{K\in\mathcal{T}_h} \int_K
       \Big|\sum_{i=0}^{N(h)} \Delta\alpha_i(\omega)\,\phi_i(x)\Big|^2\,dx \\
  &= \sum_{i,j=0}^{N(h)} \Delta\alpha_i(\omega)\,\Delta\alpha_j(\omega)
       \sum_{K\in\mathcal{T}_h} \int_K \phi_i(x)\,\phi_j(x)\,dx \\
  &= \Delta\alpha(\omega)^\top M_h\,\Delta\alpha(\omega),
\end{align*}

where the (global) DG mass matrix $M_h \in \mathbb{R}^{(N(h)+1)\times (N(h)+1)}$
is defined by
$$
  (M_h)_{ij}
  := \sum_{K\in\mathcal{T}_h} \int_K \phi_i(x)\,\phi_j(x)\,dx,
  \qquad 0 \le i,j \le N(h).
$$

Since the matrix $M_h$ is symmetric positive definite, all its eigenvalues are strictly positive. 
Denoting by $\lambda_{\max}(M_h)$ its largest eigenvalue, we obtain from 
Lemma~\ref{lem:SPD_quadratic_bound} that
$$
  \Delta\alpha(\omega)^\top M_h\,\Delta\alpha(\omega)
  \;\le\; \lambda_{\max}(M_h)\,\|\Delta\alpha(\omega)\|_2^2.
$$

Thus there exists a constant $C_h>0$, depending only on $h$ and the basis,
such that
\begin{equation}\label{eq:broken_L2_vs_coeff}
  \|\Delta u^{\mathrm{DG}}_{\varepsilon,h}(\omega,\cdot)\|_{L^2(\mathcal T_h)}^2
  \;\le\; C_h\,\|\Delta\alpha(\omega)\|_2^2
  \qquad \forall\,\omega\in\Omega_{\text{param}}.
\end{equation}
Next, we integrate over the parameter space $\Omega_{\text{param}}$.
\begin{align*}
  \| u^{\mathrm{DG}}_{\varepsilon,h}
     - \widehat{u}^{\mathrm{DG}}_{\varepsilon,h,n,M}
  \|_{L^2(\Omega_{\text{param}};L^2(\Omega))}^2
  &= \int_{\Omega_{\text{param}}}
     \| \Delta u^{\mathrm{DG}}_{\varepsilon,h}(\omega,\cdot)\|_{L^2(\mathcal{T}_h)}^2
     \,d\mathbb{P}_{\Omega_{\text{param}}}(\omega) \\
  &\le \int_{\Omega_{\text{param}}}
       C_h\,\|\Delta\alpha(\omega)\|_2^2\,d\mathbb{P}_{\Omega_{\text{param}}}(\omega) \\
  &= C_h\,\|\alpha^* - \widehat{\alpha}(n,M)\|_{L^2(\Omega_{\text{param}})}^2 \rightarrow 0.
\end{align*}
This completes the proof.
\end{proof}

\section{Numerical Experiments}{\label{sec:experiments}}


In this section, we evaluate the performance of the proposed DG-FEONet framework through several numerical experiments. We emphasize its accuracy and generalization capability in solving parametric PDEs with discontinuous data (e.g, coefficients and source terms), when the solutions can also be discontinuous. {The performance of DG-FEONet is compared against the original FEONet architecture 
(referred to as CG-FEONet) which is based on classical CG, highlighting  the advantage of the DG formulation in resolving discontinuities.}

Our implementation combines PyTorch~\cite{paszke2019pytorch} for neural network modeling with FEniCS~\cite{alnaes2015fenics} for finite element assembly, following the methodology outlined in \cite{MR4888707}. All experiments are accelerated using GPU hardware (NVIDIA A6000) to facilitate efficient training and evaluation. For all the following examples, we randomly generate 2,000 parametric input samples for the one-dimensional cases, of which 1,000 are used as training data and 1,000 as test data. For the two-dimensional cases, we generate 1,000 parametric input samples, equally split into 500 for training and 500 for testing.
The model is trained in an data-free manner using only the residuals of the weak form, without requiring paired input-output data, as discussed in the previous sections. 
Reference solutions for both training and testing are computed using the finite element solver on a sufficiently fine mesh to serve as the ground truth.

We evaluate the performance of DG-FEONet by computing the \emph{average relative discrete} $\ell^2$\emph{-error} over a set of $ N_s $ test samples. Let $ u_{i,j}^{\text{FEM}} $ and $ u_{i,j}^{\text{nn}} $ denote the FEM solution and the neural network prediction, respectively, for the $i$-th test sample $(i=1,\ldots,N_s)$ 
at the $j$-th spatial degrees of freedom, where $j = 1, \dots, M$. 
Then the relative $\ell^2$-error is defined as:
\begin{equation}
\|  u_{i,j}^{\text{nn}} - u_{i,j}^{\text{FEM}} \|_{l^2}^{\text{Rel}}:=
 \left( \frac{ \sum_{j=1}^{M} \left( u_{i,j}^{\text{nn}} - u_{i,j}^{\text{FEM}} \right)^2 }{ \sum_{j=1}^{M} \left( u_{i,j}^{\text{FEM}} \right)^2 } \right)^{1/2},
\label{eqn:relative_error}
\end{equation}
for each $i$-th sample, and the 
average relative \(\ell^2\)-error is defined as:
\begin{equation}\label{avg_rel_error}
{E}_{\text{rel}}:=
 \frac{1}{N_s} \sum_{i=1}^{N_s} \|  u_{i,j}^{\text{nn}} - u_{i,j}^{\text{FEM}} \|_{l^2}^{\text{Rel}}.    
\end{equation}

\begin{table}[]
\centering
\caption{{Summary of numerical experiments for DG-FEONet.}}
\label{tab:exp_summary}
\resizebox{\textwidth}{!}{
\begin{tabular}{lllll}
\hline
Experiment                      & Dimension & Input                                            & Output   & Source of Discontinuity \\ \hline
Exp1~(Sec.~\ref{subsec:ex_1}) & 1D        & $f(x)=k$ with $k\in[0.1,2]$                        & $u(x)$   & $\varepsilon(x)$   \\
Exp2~(Sec.~\ref{subsec:ex_2}) & 1D        & $c(x)$ with random jumps at random locations       & $u(x)$   & $\varepsilon(x), c(x)$ and $f(x)$  \\
Exp3~(Sec.~\ref{subsec:ex_3}) & 2D        & $f(x,y)$ with random sine/cosine coefficients      & $u(x,y)$ & $\varepsilon(x,y)$ \\
Exp4~(Sec.~\ref{subsec:ex_4}) & 2D        & $\varepsilon(x,y)$ with random radius and diffusion values & $u(x,y)$ & $\varepsilon(x,y)$ \\
Exp5~(Sec.~\ref{subsec:ex_5}) & 2D        & $\varepsilon(x,y)$ with two random circles, radius and diffusion values  & $u(x,y)$ & $\varepsilon(x,y)$ \\ \hline
\end{tabular}
}
\end{table}

\paragraph{Setup}
For numerical experiment 1--2, we consider a one-dimensional case with the following equation:
\begin{subequations}\label{eqn:example_1d}
\begin{align}
-\,\varepsilon(x)\,u''(x) + b(x)\,u'(x) + c(x)\,u(x) &= f(x) \text{ in } \Omega,\label{eqn:example_1d:a}\\
u(-1) = u(1) &= 0,                                   \label{eqn:example_1d:b}
\end{align}
\end{subequations}
where $\Omega:=(-1,1)$. For these experiments, the architecture is given as a fully connected neural network (FCNN) with 3 hidden layers of 32 neurons each. We use the SiLU (Sigmoid-weighted Linear Unit) activation function, and the L-BFGS optimizer (Limited-memory Broyden-Fletcher-Goldfarb-Shanno) with a batch size of 32. 

Next, for numerical experiments 3--5, we consider a two-dimensional case with the following equation:
\begin{subequations}\label{eqn:example_2d}
\begin{align}
- \nabla \cdot (\varepsilon(x,y) \nabla u(x, y)) + \mathbf{v} \cdot \nabla u(x, y) &= f(x, y) \text{ in } \Omega, \label{eqn:example_2d:a}\\
u(x,y) &= 0 \text{ on } \partial \Omega,                       \label{eqn:example_2d:b}
\end{align}
\end{subequations}
where $\Omega:=(-1,1)^2$. In these experiments, we employ a two-dimensional convolutional neural network (CNN) with 5 hidden layers, each containing 32 neurons, and SiLU activation functions. The model is trained in an data-free manner using the L-BFGS optimizer. A detailed summary of the input parametrization and the types of discontinuous coefficients considered in each numerical experiment is provided in Table~\ref{tab:exp_summary}.

\subsection{Numerical Experiment 1: Convergence Test with Discontinuous Coefficient}
\label{subsec:ex_1}
We first consider a simple one-dimensional setting to establish baseline performance and convergence behavior.  We evaluate the performance of the proposed DG-FEONet with a discontinuous coefficient and investigate the convergence behavior of the solution. 
Here, we set $c(x)=0.001$ and $b(x)=0.01$ for \eqref{eqn:example_1d:a}, but we introduce a discontinuity in the diffusion coefficient:
$$
\varepsilon(x)=
\begin{cases}
\varepsilon_{\min}, & x<0,\\
\varepsilon_{\max}=m\,\varepsilon_{\min}, & x\ge 0,
\end{cases}
\qquad
\varepsilon_{\min}=0.01,\;\; m\in\{5,10,100\}.
$$
This discontinuity introduces a gradient jump in the solution at the interface $ x = 0 $, presenting a challenge to numerical approximation by employing a neural network and operator learning~\cite{wu2022inn}. In this case, we consider $f(x):= k \in [0.1,2]$ as the input, where each function takes a constant value randomly sampled from the interval. To assess the convergence of the error ${E}_{\text{rel}}$, we obtain $u^{\text{nn}}$ on uniform meshes with $N_e = 16, 32, 64,$ and $128$ elements, where $N_e$ denotes the number of mesh elements in the domain. 

\begin{figure}[!ht]
  \centering
  \begin{minipage}{0.48\textwidth}
    \centering
    \includegraphics[width=\linewidth]{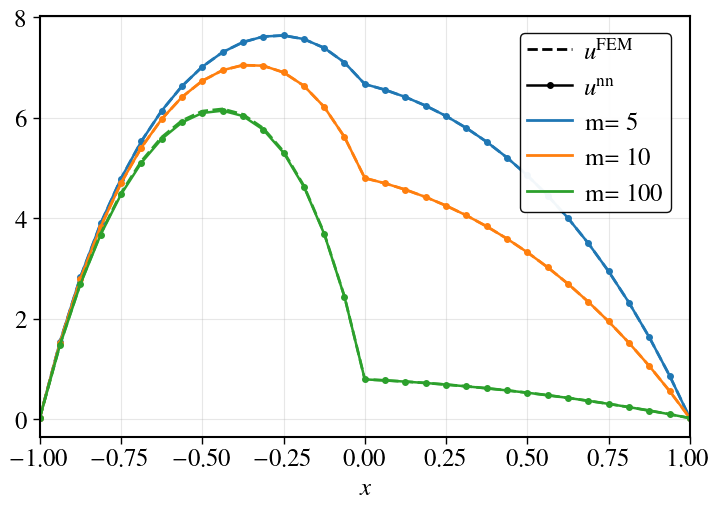}\\
    {\small (a)}
  \end{minipage}\hfill
  \begin{minipage}{0.48\textwidth}
    \centering
    \includegraphics[width=\linewidth]{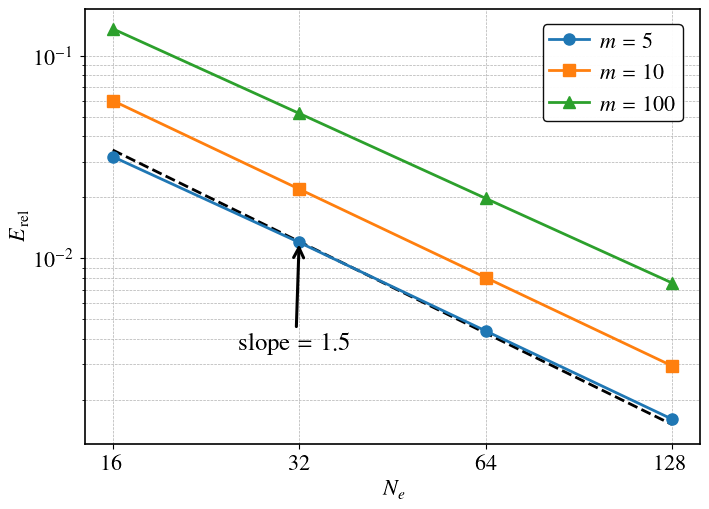}\\
    {\small (b)}
  \end{minipage}

  \caption{
  Performance of DG-FEONet for a discontinuous diffusion coefficient (Experiment~1).
  (a) Solution plots for $m=5,10,$ and $100$.
  (b) Convergence of the average relative $\ell^2$-error $E_{\mathrm{rel}}$
  for different values of $m$.
  }
  \label{fig:dis_eps_side_by_side}
\end{figure}

Figure~\ref{fig:dis_eps_side_by_side} shows both the predicted solution for a representative input and the convergence plot.
These results demonstrate the robustness of DG-FEONet in learning solutions to PDEs with discontinuous coefficients and achieving nearly sub-optimal convergence rates as expected~\cite{cai2011discontinuous}.

\begin{figure}[!ht]
    \centering
    \includegraphics[width=0.5\linewidth]{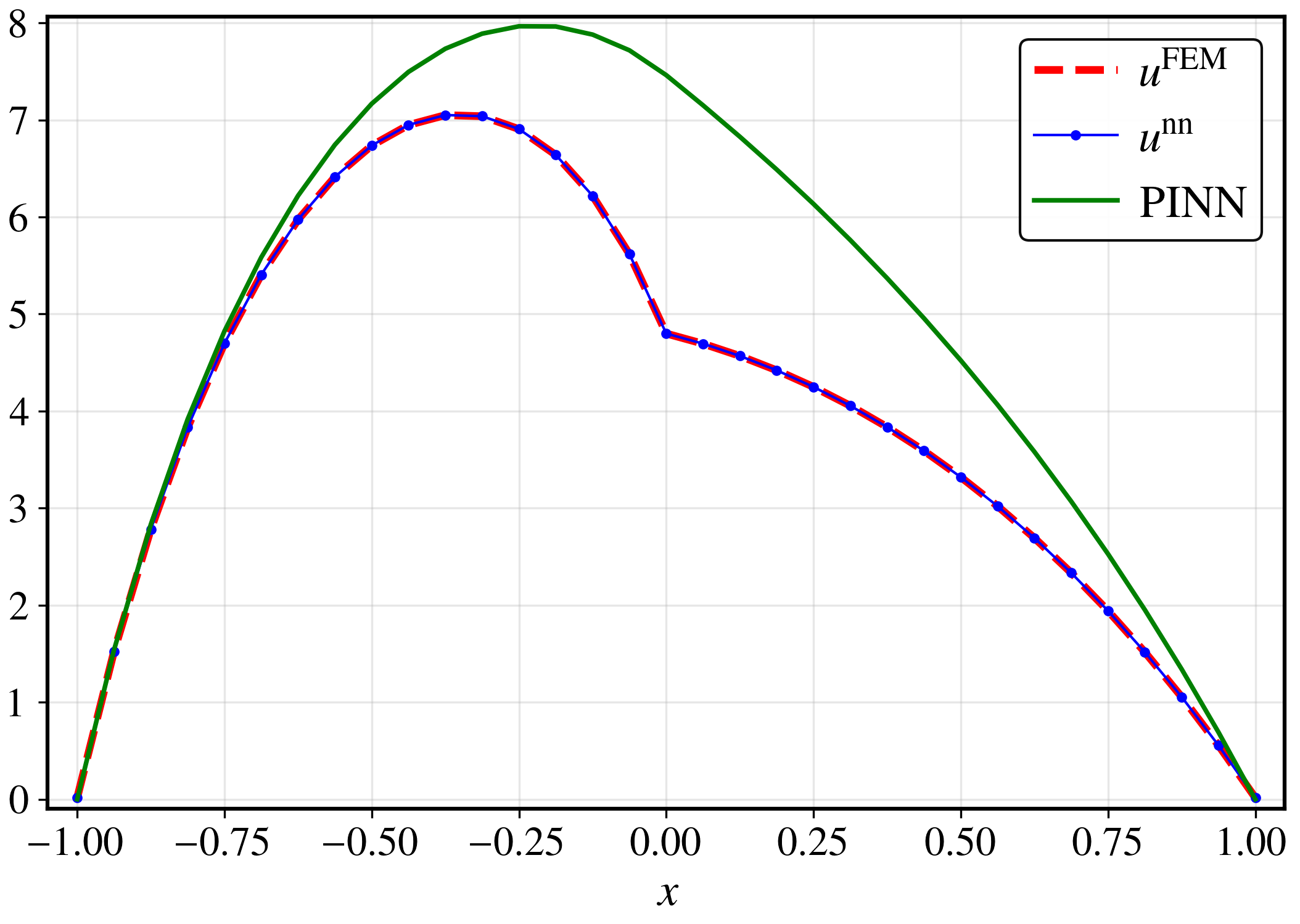}
    \caption{Comparison of the DG reference solution, DG-FEONet prediction, and PINN approximation
for Experiment~1 with a discontinuous diffusion coefficient $\varepsilon(x)$ and jump factor $m=10$.} 
    \label{fig:PINN}
\end{figure}

In addition, we include as a baseline a PINN 
with a fully connected architecture of four hidden layers. The PINN is trained on 
the same PDE~\eqref{eqn:example_1d:a} with coefficients $b(x)=0.01$, $c(x)=0.1$, 
$f(x)= 1.0$, and
$$
\varepsilon(x) =
\begin{cases}
0.01, & x < 0,\\
0.1,  & x \ge 0.
\end{cases}
$$
Figure~\ref{fig:PINN} shows that the PINN fails to capture the sharp kink 
induced by the discontinuous diffusion coefficient, whereas DG-FEONet matches 
the DG reference solution.

\subsection{Numerical Experiment 2: Randomly Located Discontinuities with Variable Jump Sizes}
\label{subsec:ex_2}
\begin{figure}[!ht]
    \centering
    \includegraphics[width=\linewidth]{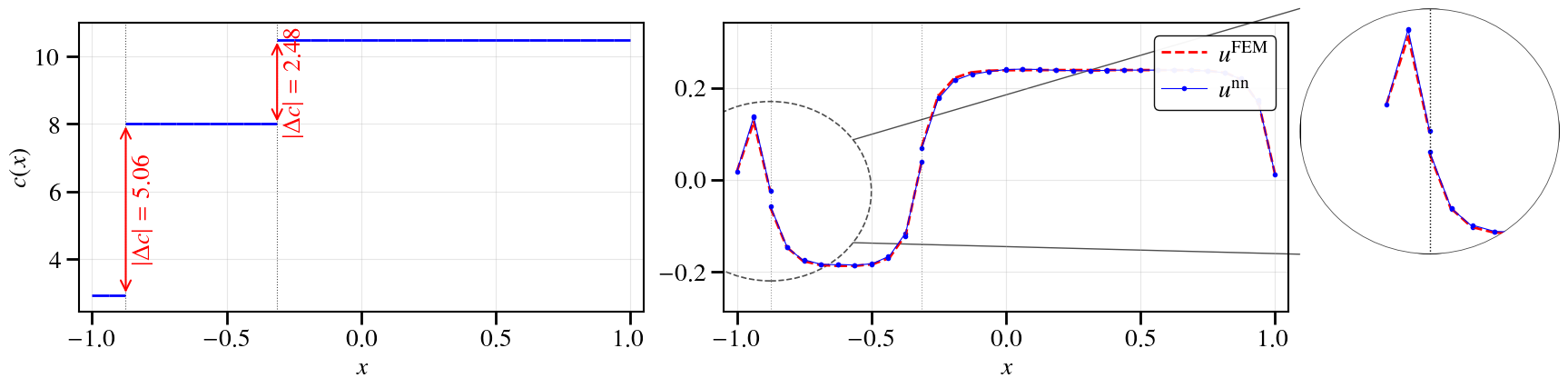}
    \includegraphics[width=\linewidth]{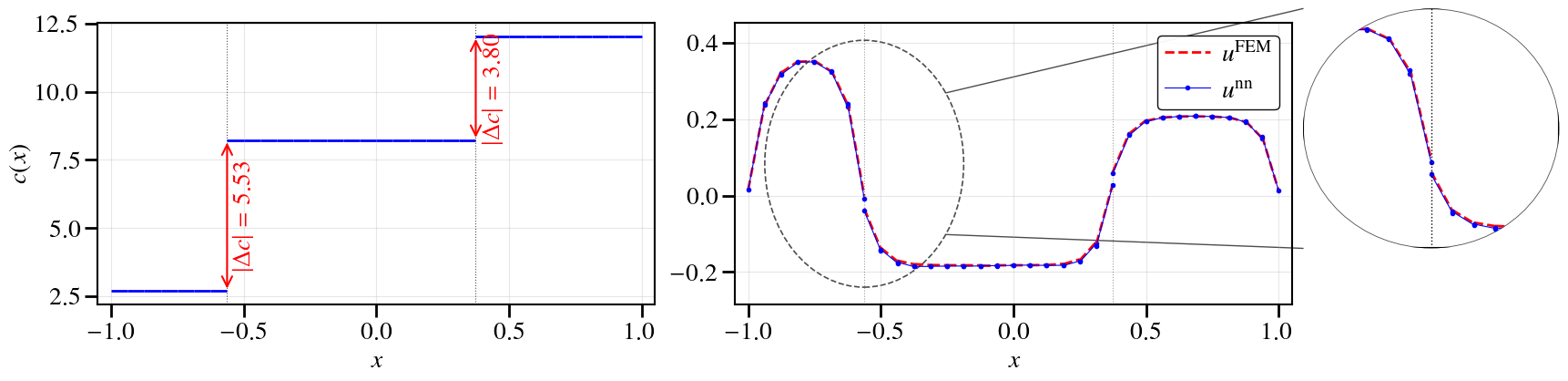}
    \includegraphics[width=\linewidth]{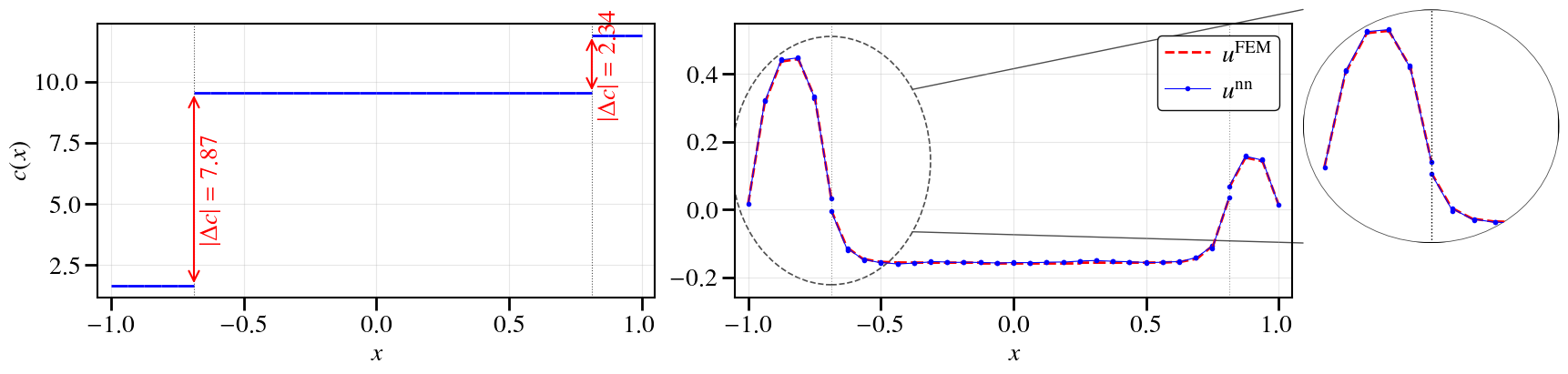}
    \caption{Three representative examples showing the input reaction coefficient $ c(x)$ with discontinuities (left) and the DG-FEONet solution $u^\text{nn}$ vs. FEM solution $ u^{\text{FEM}} $ (right) (Experiment 2). Annotated arrows indicate the magnitude of the jumps \( |\Delta c| \) across elements. DG-FEONet successfully resolves both the discontinuous input and the resulting non-smooth output.}
    \label{fig:discon_198_342}
\end{figure}

In this experiment, we consider the presence of multiple discontinuities, where both the jump magnitudes and their locations are randomly varied. In this case, the input is the reaction coefficient $c(x)$, which is defined as a random step function with two discontinuities.  
Its values are randomly sampled from three distinct intervals: $ [0,5) $, $[5,10) $, and $ [10,15] $, corresponding to the left, middle, and right subdomains, respectively.
Thus, 
$$
c(x) :=
\begin{cases}
c_0 \in [0,5), & \text{if } x < x_0, \\
c_1 \in [5,10), & \text{if } x_0 \leq x < x_1, \\
c_2 \in [10,15], & \text{if } x \geq x_1,
\end{cases}
$$
where $ x_0 $ and $ x_1 $ are also  randomly selected points in the domain $(-1,1)$. See Figure~\ref{fig:discon_198_342} for the setups of some examples.
The source term $ f(x)$ is defined as a fixed piecewise constant function given by
$$
f(x) =
\begin{cases}
1.0, & \text{if } x < x_0, \\
-1.5, & \text{if } x_0 \leq x < x_1, \\
2.5, & \text{if } x \geq x_1.
\end{cases}
$$
The diffusion coefficient is also discontinuous and is given by 
$$
\varepsilon(x,y) =
\begin{cases}
0.01, & \text{if } x < x_0, \\
0.02, & \text{if } x_0 \leq x < x_1, \\
0.03, & \text{if } x \geq x_1.
\end{cases}
$$
The other coefficient is set as $b(x)=0.01$. This setup produces diverse solution profiles, posing a challenging test for the network's generalization.

\begin{figure}[!ht]
    \centering
    \includegraphics[width=\linewidth]{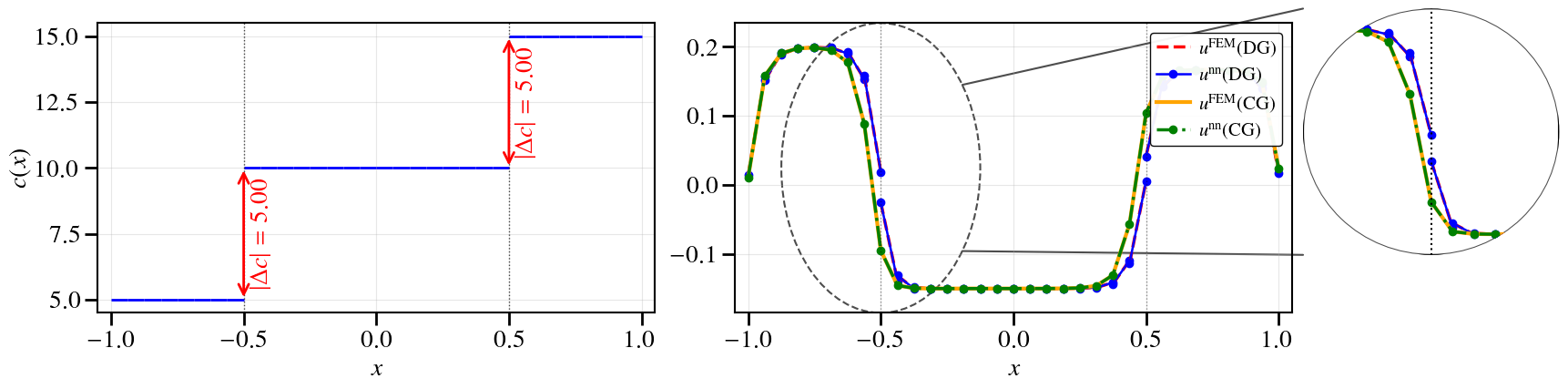}
    \caption{Comparison of CG-FEONet and DG-FEONet for a sample
    with two discontinuities (Experiment 2). The CG-FEONet solution oversmooths the jumps,
    while DG-FEONet resolves them and agrees with the DG reference solution.}
    \label{fig:cg_dg_feonet_comparison}
\end{figure}

In Figure~\ref{fig:discon_198_342}, we visualize threes representative samples of the input function $ c(x) $ and the predicted solution $ u^{\text{nn}}(x) $. The reaction coefficient $ c(x) $ exhibits clear discontinuities, and the corresponding solution $ u^{\text{nn}}(x) $ accurately captures the sharp transitions induced by these discontinuities. The DG-FEONet predictions closely matches the FEM solution, demonstrating its ability to handle discontinuous coefficients.
For comparison, we also train a CG-FEONet with the same architecture and training data, but using a continuous Galerkin discretization. Since the CG space is continuous and cannot represent true jumps, CG-FEONet smoothen out the discontinuities, whereas DG-FEONet keeps the sharp jumps and agrees with the DG reference solution; see Figure~\ref{fig:cg_dg_feonet_comparison}.

\subsection{Numerical Experiment 3: Convergence Test with a Discontinuous Diffusion Coefficient}
\label{subsec:ex_3}

We consider the two-dimensional case with the given PDE \eqref{eqn:example_2d:a} to investigate the convergence behavior of the solutions with a discontinuous diffusion coefficient
$$
\varepsilon(x,y):=
\begin{cases}
\varepsilon_{\min}, & \text{if } \sqrt{x^{2}+y^{2}}<0.5,\\[2pt]
\varepsilon_{\max}=m\,\varepsilon_{\min}, & \text{otherwise},
\end{cases}
\qquad
\varepsilon_{\min}=0.1,\;\; m\in\{5,10,20\}.
$$
 The discontinuity reduces the regularity of the solution and induces gradient jumps aligned with the interface. Here, we set $ \mathbf{v} = (-1, 0) $, and the input functions are given as
$$
f(x, y) = m_0 \sin(n_0 x + n_1 y) + m_1 \cos(n_2 x + n_3 y),
$$
where $ m_0, m_1 \in [1, 2] $ and $ n_0, n_1, n_2, n_3 \in [0, \pi] $ are randomly sampled.
To evaluate the convergence of the average of relative $\ell^2-$ error ${E}_{\text{rel}}$ defined in \eqref{avg_rel_error}, we discretize the domain using triangular meshes with 
$ N_e = 32, 128, 512$, and $2048$  elements (corresponding to the respective following number of DOFs: $96$, $384$, $1536$, and $6144$).

\begin{figure}[!ht]
  \centering

  \includegraphics[width=\textwidth]{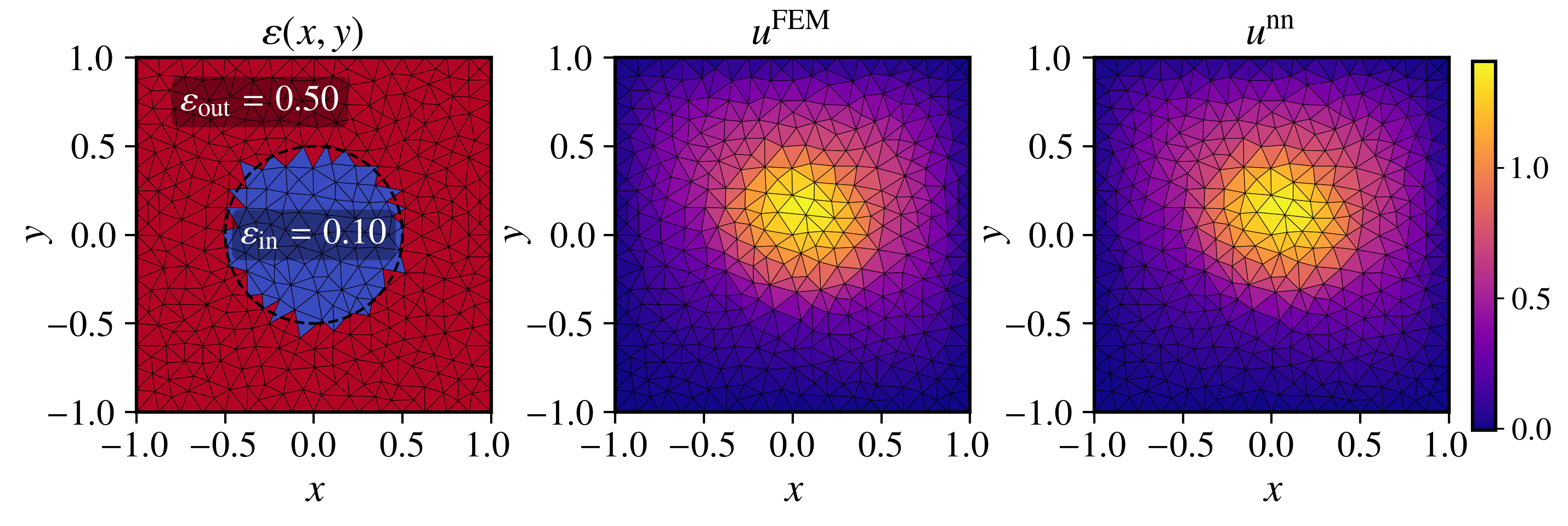}\\[-0.25em]
  {\small (a) Case 1 ($m=5$): $\varepsilon_{\mathrm{in}}=0.10$, $\varepsilon_{\mathrm{out}}=0.50$.}
  \label{fig:5x-jump}

  \vspace{0.8em}

  \includegraphics[width=\textwidth]{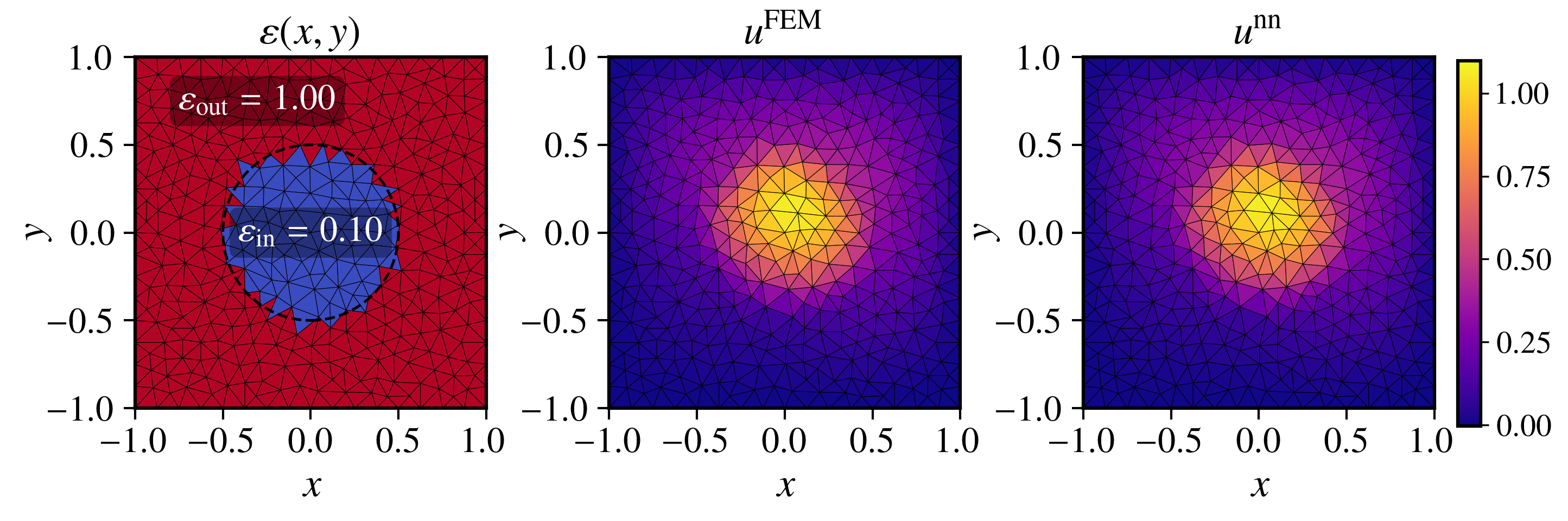}\\[-0.25em]
  {\small (b) Case 2 ($m=10$): $\varepsilon_{\mathrm{in}}=0.10$, $\varepsilon_{\mathrm{out}}=1.00$.}
  \label{fig:10x-jump}

  \vspace{0.8em}

  \includegraphics[width=\textwidth]{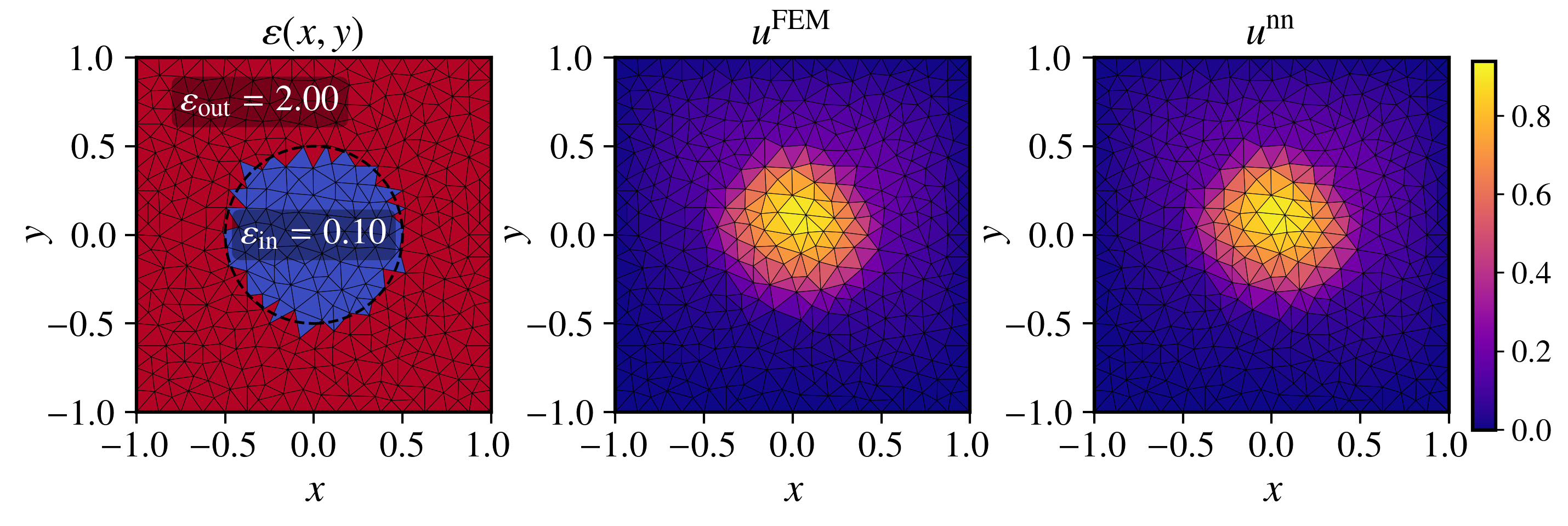}\\[-0.25em]
  {\small (c) Case 3 ($m=20$): $\varepsilon_{\mathrm{in}}=0.10$, $\varepsilon_{\mathrm{out}}=2.00$.}
  \label{fig:20x-jump}

  \caption{Element-wise plots for three diffusion jumps $\varepsilon(x,y)$ (left), FEM reference solution $u^{\mathrm{FEM}}$ (middle), and DG-FEONet prediction $u^{\mathrm{nn}}$ (right) in Experiment~3. Panels (a)–(c) correspond to $m=5$, $m=10$, and $m=20$, respectively.}
  \label{fig:dg-feonet-all}
\end{figure}

Figure~\ref{fig:dg-feonet-all} displays the output for a given representative sample for each $m\in\{5, 10, 20\}$. The first panel shows the discontinuous diffusion coefficient $ \varepsilon(x, y) $, highlighting the circular interface. The second and third panels show the FEM reference solution $ u^{\text{FEM}}(x, y) $ and the DG-FEONet prediction $ u^{\text{nn}}(x, y) $, respectively. {It confirm that DG-FEONet can stably approximate two-dimensional PDEs with discontinuous coefficients and accurately recover interface-driven solution structures across varying jump sizes.}

\begin{figure}[!ht]
  \centering
  \begin{minipage}[t]{0.48\textwidth}
    \centering
    \includegraphics[width=\linewidth]{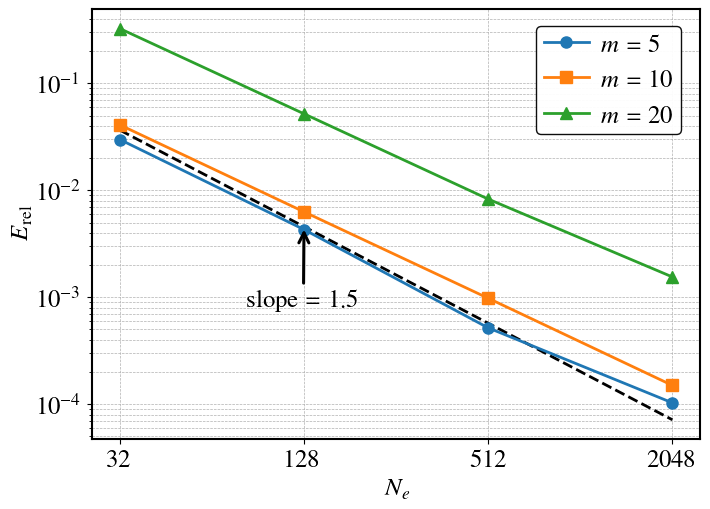}\\
    {\small (a)}
  \end{minipage}\hfill
  \begin{minipage}[t]{0.48\textwidth}
    \centering
    \includegraphics[width=\linewidth]{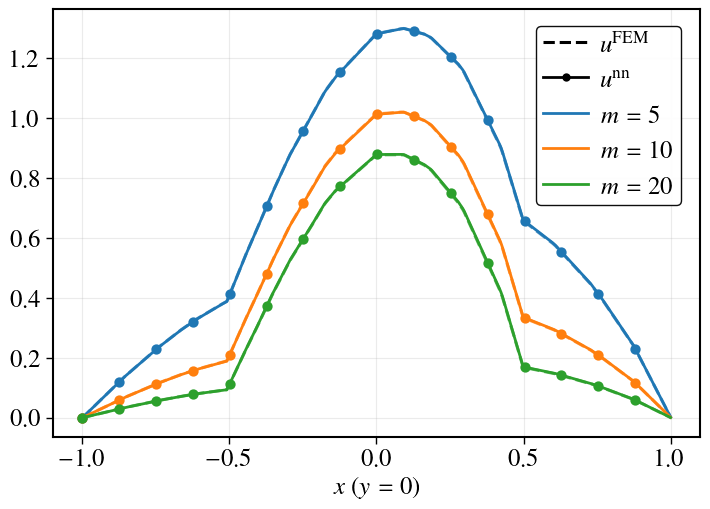}\\
    {\small (b)}
  \end{minipage}

  \caption{Convergence and solution profile for a discontinuous diffusion coefficient (Experiment~3).
  (a) Convergence of the average relative $\ell^2$-error $E_{\mathrm{rel}}$ for different values of $m$.
  (b) Line-slice comparison of $u^{\mathrm{nn}}$ and $u^{\mathrm{FEM}}$ along $y=0$ for $m=5,10,$ and $20$.}
  \label{fig:1d_discontinuity_profiles}
\end{figure}

Figure~\ref{fig:1d_discontinuity_profiles} (a) shows the expected convergence behavior under the mesh refinement with the discontinuous coefficients. The plot demonstrates that the method achieves a consistent and near-optimal convergence rate, closely following the reference slope. Figure~\ref{fig:1d_discontinuity_profiles} (b) shows the solution over the line $y=0$. We observe that the proposed DG-FEONet exhibits stable and expected convergence performance for the problems with continuous/discontinuous coefficients in 2D.

\subsection{Numerical Experiment 4: Random Circular Discontinuities with Varying Radius and Diffusion Values}\label{subsec:ex_4}






In this example, we extend the previous setting by introducing variability in both the radius and in diffusion values. The input diffusion coefficient $ \varepsilon(x, y) $ is defined as
$$
\varepsilon(x, y) = 
\begin{cases}
\varepsilon_\text{in}, & \text{if } \sqrt{x^2 + y^2} < r, \\
\varepsilon_\text{out}, & \text{otherwise},
\end{cases}
$$
where $ r \in (0, 1] $ denotes the radius, and $ \varepsilon_\text{in} \in [0.1, 0.5] $, $ \varepsilon_\text{out} \in [1.0, 2.0] $ are sampled independently for each realization. This setup generates a discontinuous diffusion coefficient with random circular interfaces and jump sizes, providing a challenging generalization test for the learned operator. Here, the source function is fixed as $f(x, y) = 1.5\sin(\pi (x + 2y)) + 1.2\cos(\pi (1.5x + 0.5y))$.

\begin{figure}[!ht]
  \centering

  \includegraphics[width=\textwidth]{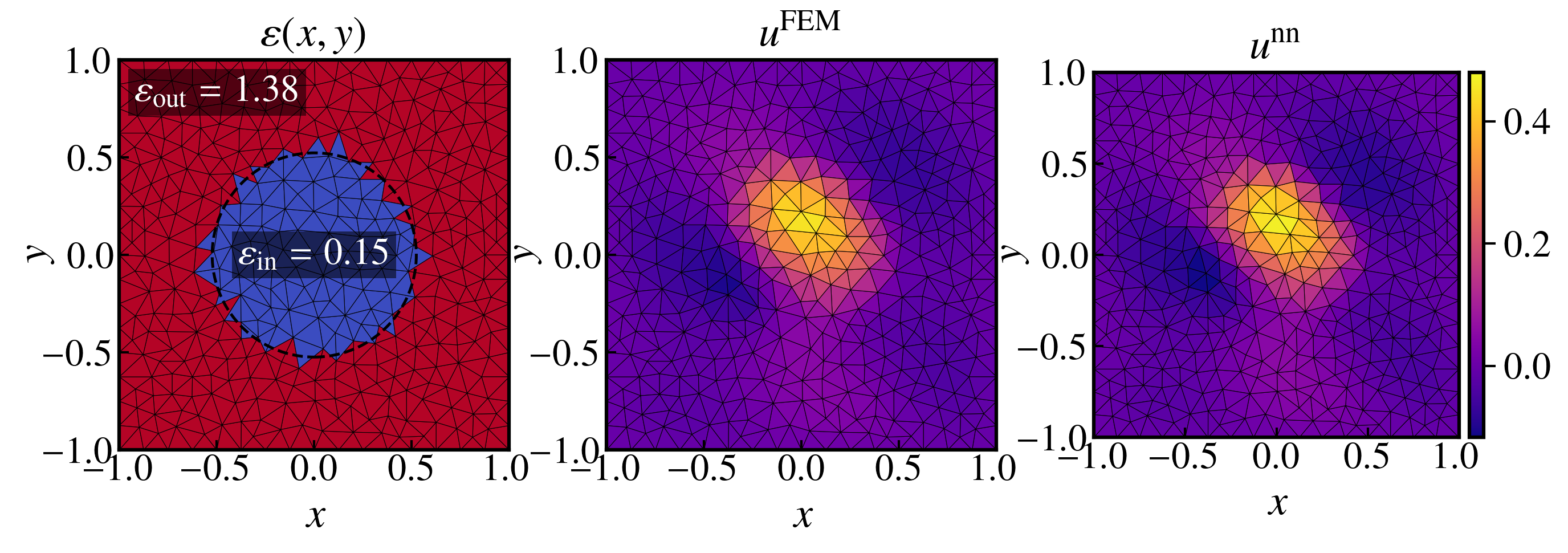}\\[-0.25em]
  {\small (a) $r=0.52$, $\varepsilon_{\mathrm{in}}=0.15$, $\varepsilon_{\mathrm{out}}=1.38$.}

  \vspace{0.9em}

  \includegraphics[width=\textwidth]{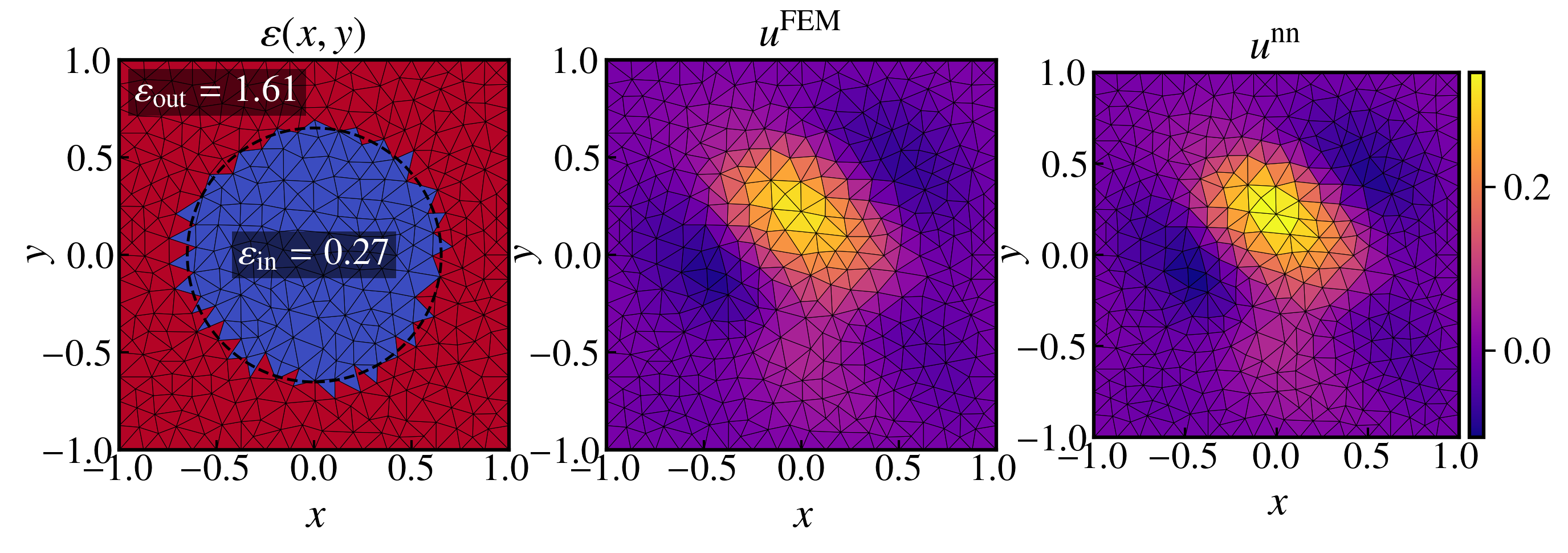}\\[-0.25em]
  {\small (b) $r=0.65$, $\varepsilon_{\mathrm{in}}=0.27$, $\varepsilon_{\mathrm{out}}=1.61$.}

  \caption{Visualization of representative diffusion realizations with random circular interfaces (Experiment~4).
  Each row shows (left) the discontinuous diffusion coefficient $\varepsilon(x,y)$, (middle) the FEM reference solution
  $u^{\mathrm{FEM}}$, and (right) the DG-FEONet prediction $u^{\mathrm{nn}}$ for the indicated
  $(r,\varepsilon_{\mathrm{in}},\varepsilon_{\mathrm{out}})$ values.}
  \label{fig:dg_feonet_three_samples}
\end{figure}

Illustrative examples of randomly generated diffusion fields are shown in Figure~\ref{fig:dg_feonet_three_samples}. Figure~\ref{fig:dg_feonet_three_samples} visualizes the model prediction for two representative test case. We display the discontinuous diffusion coefficient $ \varepsilon(x, y) $, the FEM solution $ u^{\text{FEM}}(x, y) $, and the DG-FEONet prediction $ u^{\text{nn}}(x, y) $. The results highlight the model's capacity to capture the interface-induced variations and maintain high prediction accuracy across the domain.

\begin{figure}[t]
  \centering
  \begin{minipage}[t]{0.48\textwidth}
    \centering
    \includegraphics[width=\linewidth]{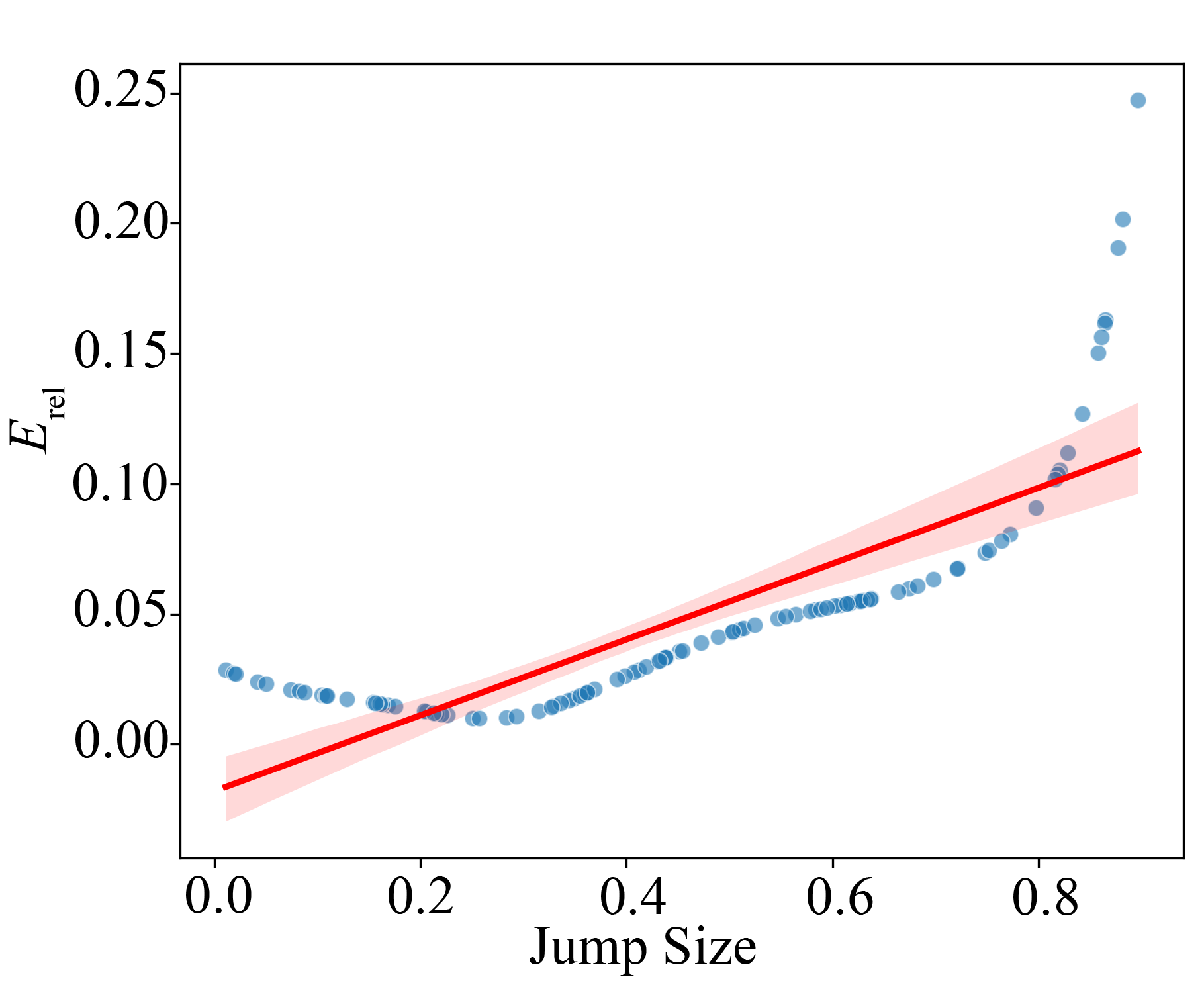}\\
    {\small (a)}
  \end{minipage}\hfill
  \begin{minipage}[t]{0.48\textwidth}
    \centering
    \includegraphics[width=\linewidth]{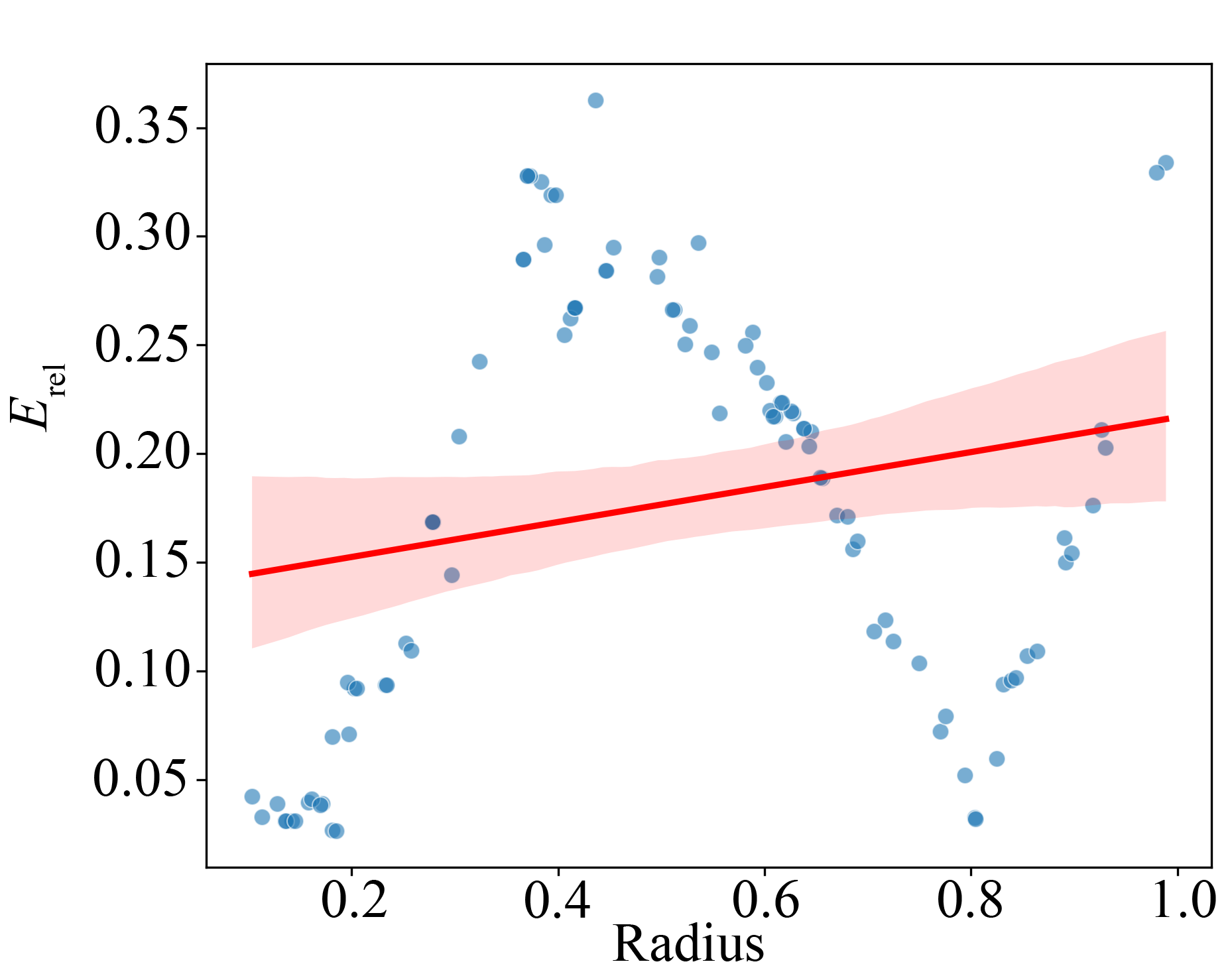}\\
    {\small (b)}
  \end{minipage}

  \caption{Sensitivity of DG-FEONet to diffusion jump size and inclusion radius (Experiment~4).
  (a) Error versus diffusion jump size $\lvert \varepsilon_{\mathrm{out}}-\varepsilon_{\mathrm{in}}\rvert$.
  (b) Error versus inclusion radius $r$.
  The red solid line denotes the least-squares regression trend, and the shaded band indicates the $95\%$ confidence interval for the fitted mean error.}
  \label{fig:jump_radius_error}
\end{figure}

To further study how the radius and jump sizes affect prediction accuracy, we conduct two experiments.  First, we fix the radius at $r=0.5$ and vary the diffusion contrast by randomly sampling $\varepsilon_{\text{in}} \in [0.1, 1.0]$, while fixing $\varepsilon_{\text{out}} = 1.0$. This isolates the effect of the diffusion jump magnitude on the model accuracy. 
Figure~\ref{fig:jump_radius_error} (a) summarizes the results, showing a clear upward trend between the relative $\ell^2$ error and the jump size 
$|\varepsilon_{\text{out}} - \varepsilon_{\text{in}}|$. 
Larger diffusion jumps lead to higher prediction errors, as expected due to stronger solution discontinuities across the interface.

Next, we fix $\varepsilon_{\text{in}} = 0.1$ and $\varepsilon_{\text{out}} = 1.0$ and vary the inclusion radius $r \in [0.1, 1.0]$ to examine how geometry (radius) influences prediction accuracy. 
Figure~\ref{fig:jump_radius_error} (b) illustrates the relationship between the relative $\ell
^2$ error and the radius. 
Although the regression line shows a slightly increasing trend, 
the scatter and wide confidence band indicate that the correlation is weak and statistically insignificant. 
Overall, DG-FEONet remains largely insensitive to variations in the inclusion radius within this range.

\subsection{Numerical Experiment 5: Two Random Non-Overlapping Circular Discontinuities with Varying Radii and Diffusion Values}
\label{subsec:ex_5}

\begin{figure}[t]
  \centering
  \includegraphics[width=\textwidth]{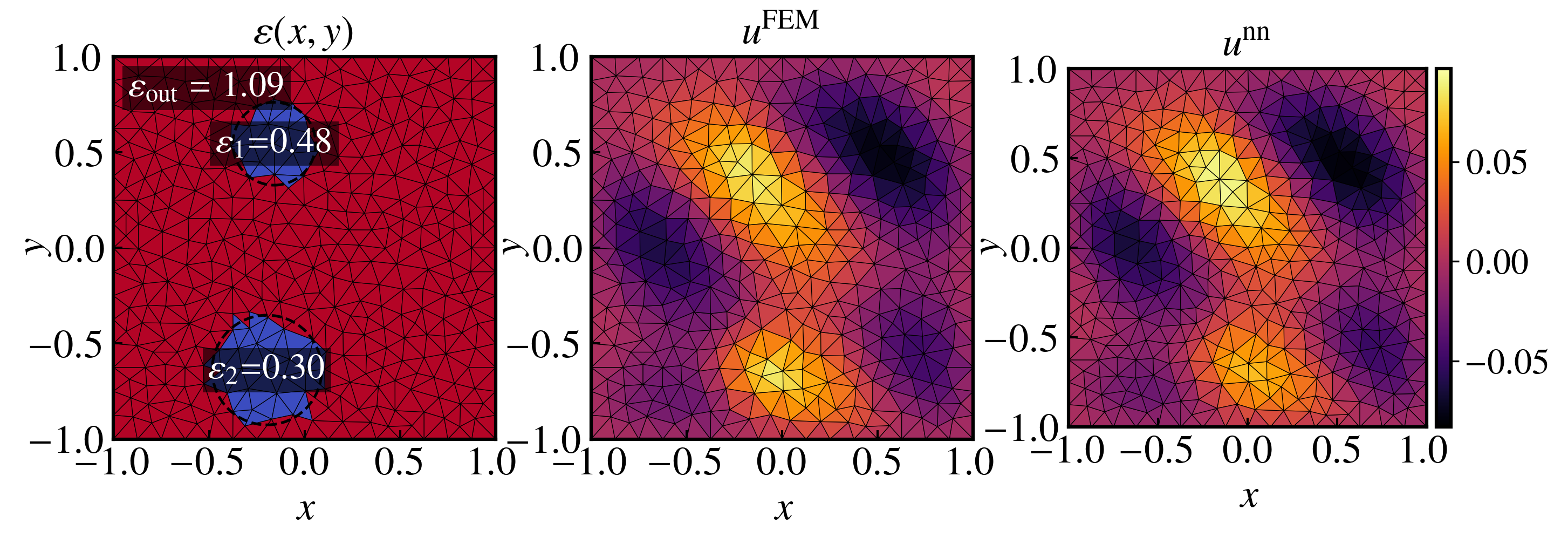}
  \caption{Representative test case with two circular inclusions of radii
  $r_1=0.169$ and $r_2=0.397$, and diffusion coefficients
  $\varepsilon_1=0.48$, $\varepsilon_2=0.30$, and $\varepsilon_{\mathrm{out}}=1.09$ (Experiment~5).
  From left to right: discontinuous diffusion field $\varepsilon(x,y)$, FEM solution $u^{\mathrm{FEM}}(x,y)$,
  and DG-FEONet prediction $u^{\mathrm{nn}}(x,y)$.}
  \label{fig:rand_radius_example_8}
\end{figure}

This example introduces two disjoint circular inclusions whose centers, radii, and diffusion values are independently randomized. The centers $(x_1, y_1)$ and $(x_2, y_2)$ are randomly sampled within the domain $\Omega = [-1, 1]^2$, subject to the constraints that the circles do not overlap and remain fully inside the domain.
Then, the input function $ \varepsilon(x, y) $ is defined as
\begin{equation}
\varepsilon(x, y) =
\begin{cases}
\varepsilon_1, & \text{if } \left|(x, y) - (x_1, y_1)\right| < r_1, \\
\varepsilon_2, & \text{if } \left|(x, y) - (x_2, y_2)\right| < r_2 \text{ and not in circle 1}, \\
\varepsilon_\text{out}, & \text{otherwise},
\end{cases}
\end{equation}
where the radii $ r_1, r_2 \in [0.2, 0.4] $, interior diffusion values $ \varepsilon_1, \varepsilon_2 \in [0.1, 0.5] $, and background diffusion $ \varepsilon_\text{out} \in [1.0, 2.0] $ are all drawn independently at random from their respective intervals. Here, we set $ \mathbf{v} = (-1, 0) $, and fixed the source term to be 
$
f(x, y) = 1.5\sin(\pi (x + 2y)) + 1.2\cos(\pi (1.5x + 0.5y))$.
Figure~\ref{fig:rand_radius_example_8} illustrates a representative test sample, showing the discontinuous diffusion coefficient 
$\varepsilon(x, y)$, the FEM reference solution $u^{\mathrm{FEM}}(x, y)$, and the DG-FEONet prediction $u^{\mathrm{nn}}(x, y)$. 
The displayed configuration includes two circular inclusions with radii 
$r_1 = 0.169$ and $r_2 = 0.397$, and diffusion coefficients 
$\varepsilon_1 = 0.48$, $\varepsilon_2 = 0.30$, and $\varepsilon_\mathrm{out} = 1.09$. 

\begin{figure}[t]
  \centering
  \includegraphics[width=\textwidth]{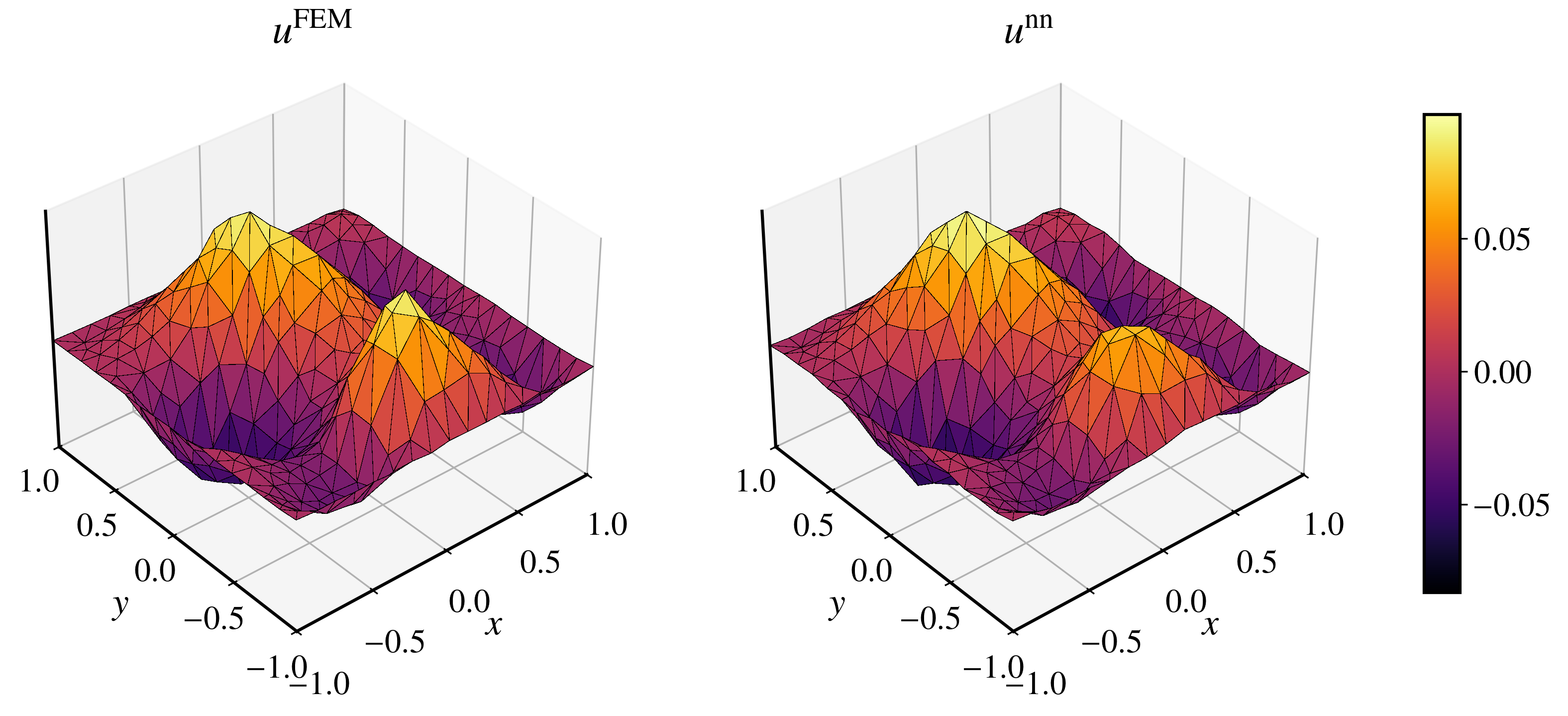}
  \caption{Three-dimensional surface comparison for Experiment~5.
  Left: FEM reference solution $u^{\mathrm{FEM}}(x,y)$.
  Right: DG-FEONet prediction $u^{\mathrm{nn}}(x,y)$.}
  \label{fig:rand_radius_example_8_3D}
\end{figure}

Figure~\ref{fig:rand_radius_example_8_3D} shows the 3D surface comparison of 
$u^{\mathrm{FEM}}(x, y)$ and $u^{\mathrm{nn}}(x, y)$, 
demonstrating close agreement in both amplitude and spatial structure, 
with only minor local smoothing observed in the DG-FEONet prediction for this complex setup.

\section{Conclusion}{\label{sec:conclusion}}

In this work, we presented DG-FEONet, a physics-informed neural operator framework based on the DG method for solving parametric PDEs. By embedding the weak formulation of the DG discretization into the loss function, DG-FEONet enables data-free learning directly from the governing equations, eliminating the need for paired input-output training data.

We evaluated the performance of DG-FEONet on a range of one- and two-dimensional advection–diffusion problems involving smooth and discontinuous coefficients. The numerical results demonstrate that DG-FEONet accurately captures sharp solution features induced by coefficient jumps, even under significant variation in interface geometry, location, and magnitude. The model exhibits favorable convergence behavior and generalization capabilities, particularly in complex settings with multiple discontinuities.

Our experiments further confirm that the DG-based architecture is especially effective in handling non-smooth problems, where traditional neural solvers and collocation based physics-informed methods often struggle. By leveraging the local nature and flexibility of the DG method, DG-FEONet provides a scalable and robust approach to PDE learning in irregular and heterogeneous environments.

\appendix

\section{Invertibility of the DG operator matrix}\label{app:A_invertible}

\begin{lemma}[Invertibility of the DG operator matrix]\label{lem:A_invertible}
Assume that the DG bilinear form $a(\cdot,\cdot)$ is coercive in the DG energy norm,
i.e., there exists $C_{\mathrm{coer}}>0$ such that
$$
a(v_h,v_h)\ge C_{\mathrm{coer}}\|v_h\|_{DG}^2
\quad \forall v_h\in V_h.
$$
Let $A\in\mathbb{R}^{(N_h+1)\times(N_h+1)}$ be defined by
$$
A_{ij}:=a(\phi_j,\phi_i).
$$
Then $A$ is non-singular.
\end{lemma}

\begin{proof}
    Let $\alpha = (\alpha_0,\dots,\alpha_{N_h})^\top \in \mathbb{R}^{N_h+1}$ satisfy
$$
A\alpha = 0.
$$
Define the corresponding DG function $u_h \in V_h$ by
$$
u_h := \sum_{j=0}^{N_h} \alpha_j \,\phi_j.
$$
By the definition \eqref{eq:A_def_main} of $A$, the $i$-th component of $A\alpha$ is
$$
(A\alpha)_i = \sum_{j=0}^{N_h} A_{ij}\,\alpha_j
            = \sum_{j=0}^{N_h} a(\phi_j,\phi_i)\,\alpha_j
            = a\Big(\sum_{j=0}^{N_h} \alpha_j \phi_j,\ \phi_i\Big)
            = a(u_h,\phi_i),
$$
for all $i=0,\dots,N_h$. Thus $A\alpha = 0$ implies
$$
a(u_h,\phi_i) = 0 \qquad \forall\, i=0,\dots,N_h.
$$
Since $\{\phi_i\}_{i=0}^{N_h}$ is a basis of $V_h$, any $v_h \in V_h$ can be written as
$$
v_h = \sum_{i=0}^{N_h} \beta_i \phi_i.
$$
By bilinearity of $a(\cdot,\cdot)$, we have
$$
a(u_h,v_h) = \sum_{i=0}^{N_h} \beta_i\,a(u_h,\phi_i) = 0
\qquad \forall\, v_h \in V_h.
$$
In particular, taking $v_h = u_h$ yields
$$
a(u_h,u_h) = 0.
$$
Using the coercivity assumption \eqref{eq:coercivity_main}, we obtain
$$
0 = a(u_h,u_h) \;\ge\; C_{\mathrm{coer}} \,\|u_h\|_{DG}^2.
$$
Since $C_{\mathrm{coer}}>0$, it follows that
$$
\|u_h\|_{DG}^2 = 0,
$$
and thus $u_h \equiv 0$ in $\Omega$ because $\|\cdot\|_{DG}$ is a norm on $V_h$.

Finally, $u_h \equiv 0$ and the expansion $u_h = \sum_{j=0}^{N_h} \alpha_j \phi_j$ together with the linear independence of $\{\phi_j\}_{j=0}^{N_h}$ imply
$$
\alpha_j = 0 \qquad \forall\, j=0,\dots,N_h,
$$
i.e., $\alpha = 0$. Therefore, the only vector in the kernel of $A$ is the zero vector, so $A$ is non-singular and hence invertible.

\end{proof}

\section*{Acknowledgments}
The work of K. Chawla and S. Lee was supported by the U.S. Department of Energy, Office of Science, Energy Earthshots Initiatives under Award Number DE-SC 0024703.
The work of Y. Hong was supported by Basic Science Research Program through the National Research Foundation of Korea (NRF) funded by the Korea government~(MSIT) (RS-2023-00219980), and by Institute of Information \& communications Technology Planning \& Evaluation (IITP) grant funded by the Korea government~(MSIT) [NO.RS-2021-II211343, Artificial Intelligence Graduate School Program (Seoul National University)].
The work of J.Y. Lee was supported by the Institute of Information \& Communications Technology Planning \& Evaluation (IITP) grant funded by the Korea government (MSIT) [RS-2021-II211341, Artificial Intelligence Graduate School Program (Chung-Ang University)] and by the National Research Foundation of Korea(NRF) grant funded by the Korea government(MSIT) (RS-2025-02303239).

\bibliographystyle{siamplain}
\bibliography{references}  

\end{document}